\documentclass[11pt]{amsart}

\usepackage{hyperref}
\usepackage{diagbox}
\usepackage{ upgreek }
\usepackage{ulem}

\usepackage{enumerate}
\newcommand{\hh}{{\hspace{.3mm}}}
\newcommand{\mm}{{\hspace{-.3mm}}}

\usepackage{amsmath}
\usepackage{lmodern}

\usepackage[T1]{fontenc}
\usepackage[latin9]{inputenc}
\usepackage{subfigure,lmodern, amsmath,amssymb, graphicx, pifont, adjustbox, bm, xcolor}
\usepackage{amsfonts}

\usepackage{comment}
\usepackage{mathtools}
\usepackage{float}
\usepackage{slashed}
\usepackage{ragged2e}
\usepackage{bm}
\usepackage{mathrsfs}

\newcommand{\cC}{\boldsymbol{C}}
\newcommand{\cg}{\boldsymbol{g}}
\newcommand{\csigma}{\boldsymbol{\sigma}}
\newcommand{\ceta}{\boldsymbol{\eta}}
\newcommand{\ctau}{\boldsymbol{\tau}}

\newcommand{\cA}{\boldsymbol{A}}
\newcommand{\cB}{\boldsymbol{B}}

\DeclareMathOperator{\EXT}{d}
\newcommand{\ext}{{\EXT\hspace{.01mm}}}

\newcommand{\sss}{\scriptscriptstyle}
\newcommand{\pdot}{{\textstyle\boldsymbol \cdot}\hspace{.05mm}}
\newcommand{\nn}[1]{(\ref{#1})}

\newtheorem{theorem}{Theorem}[section]

\theoremstyle{definition}

\theoremstyle{remark}
\newtheorem{remark}[theorem]{Remark}

\newtheorem{problem}[theorem]{Problem}

\newcommand{\II}{{\rm {  I\hspace{-.2mm}I}}}
\newcommand{\IIo}{\hspace{.1mm}{\rm \mathring{{ I\hspace{-.2mm} I}}}{\hspace{.2mm}}}
\newcommand{\IIIo}{\hspace{.7mm}{\rm\mathring{\!\hh\hh{ I\hspace{-.2mm} I\hspace{-.2mm}I}}{\hspace{.2mm}}}}
\newcommand{\IVo}{{
\hspace{.6mm}
\rm {\mathring{\!{ I\hspace{-.7mm} V}}}\hspace{-.3mm}}}

\newcommand{\cIIo}{\hspace{.1mm}{ \mathring{{ \boldsymbol{I}\hspace{-.7mm}  \boldsymbol{I}\hh\hh}}}{\hspace{-.8mm}}}
\newcommand{\cIIIo}{\hspace{.7mm}{\mathring{\!\hh\hh\hh\!{  \boldsymbol{I}\hspace{-.9mm}  \boldsymbol{I}\hspace{-.9mm} \boldsymbol{I}}}{\hspace{-.5mm}}}}
\newcommand{\cIVo}{{
\hspace{.6mm}
 {\mathring{\!{  \boldsymbol{I}\hspace{-.7mm} \boldsymbol{V}}}}\hspace{-.3mm}}}

\newcounter{mnotecount}

\newcommand{\mnotex}[1]
{\protect{\stepcounter{mnotecount}}$^{\mbox{\footnotesize $\bullet$\themnotecount}}$ 
\marginpar{
\raggedright\tiny\em
$\!\!\!\!\!\!\,\bullet$\themnotecount: #1} }

\begin{document}
\subjclass[2010]{53A30, 53A55, 53B25, 53C40, 53C21, 53B50, 35C20, 35G30, 35J25}

\renewcommand{\today}{}
\title{
 Boundary Curvature Scalars on\\Conformally Compact Manifolds
}
\author{ A. R. Gover${}^\flat$,
J. Kopi\'nski ${}^\natural$
 \&  A Waldron${}^\sharp$}

\address{${}^\flat$
  Department of Mathematics\\
  The University of Auckland\\
  Private Bag 92019\\
  Auckland 1142\\
  New Zealand
} \email{gover@math.auckland.ac.nz}

  \address{${}^{\natural,\sharp}$
  Center for Quantum Mathematics and Physics (QMAP)\\
  Department of Mathematics\\ 
  University of California\\
  Davis, CA95616, USA} \email{jkop,wally@math.ucdavis.edu}

\vspace{10pt}

\renewcommand{\arraystretch}{1}

\begin{abstract}
We introduce a sequence of conformally invariant scalar curvature quantities,
defined along the conformal infinity of a conformally compact~(CC) manifold, that measure the failure of a CC metric to have
constant negative scalar curvature in the interior, {\it i.e.} its failure to solve the singular Yamabe problem. Indeed, these ``CC  boundary curvature scalars'' compute  canonical expansion coefficients for  singular Yamabe metrics. Residues of their poles  
yield obstructions to smooth solutions to the singular Yamabe problem and thus, in particular, give an alternate derivation of generalized Willmore invariants. Moreover, in a given dimension, the critical CC boundary scalar characterizes the image of a Dirichlet-to-Neumann map for the singular Yamabe problem. We give explicit formul\ae\ for the first five~CC  boundary curvature scalars required for a global study of four dimensional  singular Yamabe metrics, as well as asymptotically de Sitter spacetimes.

%
%
%
%
%
%

\vspace{10cm}

\noindent
\begin{keywords}
{}\sf \tiny Keywords: Singular Yamabe problem, conformal geometry, conformally compact, 
Dirichlet-to-Neumann, 
boundary and hypersurface invariants, 
Willmore type energies, asymptotically de Sitter,
AdS/CFT, general relativity.
\end{keywords}

\end{abstract}

\maketitle

\pagestyle{myheadings} \markboth{Gover, Kopi\'nski \& Waldron}{Conformally compact boundary curvatures}

\newpage



\section{Introduction}

Especially when causal or conformal structures are of chief interest, when possible
it is propitious to study complete metrics on non-compact manifolds by
mapping this data diffeomorphically to the interior of
a compact manifold with boundary.  The image of the metric under this
map does not extend in a non-singular way to the boundary, so one
seeks to understand and exploit aspects of the structure that do
extend. A key example of such a singular metric is that of a {\it conformally
  compact structure}~$(M_+,g_+)$, where $M_+$ is the interior of a
compact, for simplicity oriented, smooth $d$-manifold $M^d$, with $\bar
d:=d-1$ dimensional boundary $\Sigma^{\bar d}:=\partial M$, and the
``singular metric'' $g_+$ is {\it conformally compactifiable} meaning
that given any smooth defining function $\eta$ for $\Sigma$, the
corresponding ``compactified metric''

\begin{equation}\label{compactified metric}
g=\eta^2 g_+
\end{equation}
smoothly extends to a metric on $M=\overline {M_+}$. Recall here that a 
smooth defining function $\eta\in C^\infty M$ obeys $\eta^{-1}(0)=\Sigma$ with nowhere vanishing conormal $\ext \eta|_\Sigma$.  We will  focus on the case where $g$ is a Riemannian metric. All our tensorial results can be easily translated to other metric signatures, but global aspects can differ markedly.

Given any strictly positive function $0<\Omega\in C^\infty M$, then $(\Omega \eta)^2 g_+$ is also a compactified metric, so we may rewrite Equation~\nn{compactified metric} as
$$
\cg  = {\ceta} ^2g_+\, ,
$$
where $\cg=[g]$ denotes a conformal class of metrics on $M$ and $\ceta = [g,\eta]$ is a {\it defining  density}. 
Recall that such a {\it conformal density} 
 may either be viewed as an equivalence class of metric-function pairs $(g,\eta)\sim (\Omega^2 g, \Omega \eta)$ or as a section of the $(1/d)^{\rm th}$ power of the bundle of volume forms tensored by the trivial line bundle,  denoted~${\mathcal E}M[1]$; see~\cite{Curry} for a review of key conformal geometric constructions. 
 Replacing~$1/d$ by~$w/d$,  the corresponding {\it weight $w$ conformal density bundle} is denoted by ${\mathcal E}M[w]$. Note that $\Gamma({\mathcal E}M[0])=C^\infty M$.
A conformal density~$\ceta$ is {\it defining} when any representative~$\eta$  of~$\ceta$ is a defining function for $\Sigma$. If any representative function~$\tau$ for  $\ctau \in \Gamma( {\mathcal E}M[1])$ is non-negative, we call the density $\ctau$ a {\it scale}. Note that each scale $\ctau$ determines a representative metric $g=\ctau^{-2}\cg$ of $\cg$.
Given a weight-$w$
 conformal-density-valued tensor $\boldsymbol{T}$, we will often write the corresponding unbolded symbol $T$ to denote any representative in a given choice of metric $g\in \cg$. 
 When 
 $T\in \Gamma({\mathbb T} M)$, for some vector bundle ${\mathbb T} M$, then we denote the corresponding weighted bundle ${\mathbb T}M[w]:={\mathbb T}M\otimes {\mathcal E}M[w]$.
Also, if further disambiguation for some equation is required, we will place a $g$ over the equality symbol   $\stackrel{g}=$, to indicate which metric $g\in \cg$ was used.

\smallskip

A key property of smooth, Riemannian-signature, CC manifolds~\cite{Loewner,Aviles,ACF,Maz} is that they admit a unique,  regularity (at least) $C^d M$,  defining density $\csigma$ that is
smooth on $M_+$ and such that the scalar curvature of the corresponding singular metric on $M_+$ has constant negative scalar curvature
 $$
Sc^{\csigma^{-2}\cg}=-d(d-1)\, .
$$
The metric $g_+^{\rm SY}:=\csigma^{-2}\cg$ is termed the {\it singular Yamabe metric}.
The {\it singular Yamabe problem} amounts to finding the {\it singular Yamabe defining density}~$\csigma$  
given the data $(M,\cg)$. We may view   the conformally embedded hypersurface $\Sigma\hookrightarrow (M,\cg)$ itself, where $\Sigma=\partial M$, as  ``Dirichlet data'' for the singular Yamabe problem.

It is well known that the scalar curvature of a CC metric extends smoothly to the boundary.  When this additionally
 obeys
$$Sc^{g_+^{\sss\rm AH}}\big|_\Sigma = -d(d-1)\, ,$$
where we have slightly abused notations,
 then $g_+^{\rm AH}$ it is termed {\it asymptotically hyperbolic}.
A singular Yamabe metric $g_+^{\rm SY}$  subject to the Einstein condition $$P^{g_+^{\rm SY}}=-\frac12 g_+^{\rm SY}\,, $$  is termed a 
 {\it Poincar\'e--Einstein metric}. Here $P$ denotes the Schouten tensor and our tensor conventions are given in Section~\ref{conv}. 
 
 \smallskip
 
There are now two natural questions we can ask. Firstly, since the singular Yamabe metric is completely determined by the embedding data $\Sigma\hookrightarrow (M,\cg)$, we can try to characterize the latter. Locally this classification question is, in part, answered by the conformal fundamental forms of~\cite{BGW}. These  arise upon asking whether the metric~$g_+^{\rm SY}$ is conformally related to a Poincar\'e--Einstein one;
see Section~\ref{CFFs}.
The second question, and the focus of this article, is to conformally invariantly characterize the failure of a CC metric $g_+$ to be the singular Yamabe metric $g_+^{\rm SY}$.

Working in analogy with the conformal fundamental form construction, one might think that this problem ought be handled by studying  jets, transverse to the boundary, of the extension to $\Sigma$ of the scalar curvature $Sc^{g_+}$. 
However, we would like a uniform construction that also handles  a second problem that cannot be dealt with this way. Namely there is a notion of a Dirchlet-to-Neumann map for singular Yamabe solutions which we shall now discuss. For that, as mentioned above,  one first views the zero locus $\Sigma$ of $\csigma$ as Dirichlet data. Then, given any $C^\infty$ smooth  defining density $\csigma_0$
for $\Sigma$, 
one can find smooth densities 
$\cA_i\in \Gamma({\mathcal E}M[-i])$ 
for 
$i\in \{0,\ldots ,\bar d\}$ such that
$
\csigma_{\rm as} = \csigma_0 \sum_{i=0}^{\bar d} \cA_i \csigma_0^i
$
obeys~\cite{ACF,WillB,CAG,Will2}
$$
Sc^{\csigma_{\rm as}^{-2}\cg}=-d(d-1)\big(1 +  \csigma_{\rm as}^d \,\cB \big)\, ,
$$
for some smooth density $\cB \in \Gamma({\mathcal E}M[-d])$. Moreover, the singular Yamabe defining density
\begin{equation}\label{C0}
\csigma = \csigma_{\rm as} + \tfrac{d}{2(d+1)}\cB\hh  \csigma_{\rm as}^{d+1}\log(\csigma_{\rm as} /\ctau) +  \csigma_{\rm as}^{d+1}\hh\cC+{\rm hots}\, ,
\end{equation}
where $\ctau$ is some scale,  ${\rm hots}$ denotes terms of higher order in  a polyhomogeneous expansion in powers of $\csigma_{\rm as}$ and $\log(\csigma_{\rm as} /\ctau)$,
and the weight $-d$ density ${\bm C}$  captures  Neumann data but is not  canonical. 
Also from the above, we see that $\cB$ is an obstruction to $C^\infty$ solutions to the singular Yamabe operator, so in dimension~$d$, one defines
the {\it obstruction density}
$$
\bar \cB_d:=\cB|_\Sigma \in \Gamma({\mathcal E}\Sigma[-d])\, . 
$$
In particular, see respectively~\cite{ACF,CAG} and~\cite{hypersurface_old}, 
\begin{align}\label{B3}
\bar B_{3} &\stackrel g= -\tfrac13 {\sf L}^{(3)}_{ab} \IIo^{ab}\,  ,\\[1mm]
\label{B4}
\bar B_{4} &\stackrel g= \tfrac16\Big(
\bar {\sf L}^{ab} \big[\IIIo_{ab}
+2 \IIo_{c(a}\IIo_{b)\circ\!}{}^c\big]
-\hh \IIo^{ab}\IVo{}_{ab}^{(4)}
+\IIIo^{ab} \IIo_{ac}\IIo^c{}_b
+2\IIIo^{ab}\IIIo_{ab}
\\&\qquad\:\:
+\tfrac12K^2
+ W^{\hat n abc} (W_{\hat n abc})^\top
\Big)\, ,\nonumber
\end{align}
where the conformally invariant second order differential operators $ {\sf L}^{(3)}_{ab}$ and $\bar {\sf L}^{ab}$ are defined in and under Equation~\nn{Lop} while our extrinsic geometry notations are given in Section~\nn{CFFs}.
Note that $\bar \cB_3$ is termed the {\it Willmore invariant} since it is the gradient of the
Willmore surface functional, with respect to variations of embeddings. Note that $\bar \cB_4$ 
and indeed $\bar \cB_d$ are  also variational, see~\cite{hypersurface_old} and~\cite{GrVol,RenVol}.
In Section~\ref{DN} we will give a canonical Dirichlet-to-Neumann map.
This will be constructed from linear operators that act on defining densities.

\smallskip

We are thus lead to the 
 following problem. 
\begin{problem}\label{prob}
 Let $(M,\cg,\ceta)$ be a CC $d$-manifold and $k\in{\mathbb Z}_{>0}$. Find conformally invariant normal operators~$\delta_k$ defined on  the space of weight one densities  $\Gamma({\mathcal E}_0 M[1])\ni {\boldsymbol f}$ that vanish along~$\Sigma$ such that they satisfy: 
 \begin{enumerate}[(i)]
\item\label{NO} {\it Normal order $k$}: 
$$\delta_k f \stackrel g= \hat n^{a_1}\cdots \hat n^{a_k} \nabla_{a_1}\cdots \nabla_{a_k} f|_\Sigma+{\rm ltots}\, ,
$$
where ``ltots'' denotes terms of lower transverse derivative order acting on $f$.\\
\item \label{CW} {\it Conformal weight:} $\Gamma({\mathcal E}_0 M[1])\stackrel{\delta_k}\longrightarrow \Gamma({\mathcal E} M[1-k ])$
so that
$\delta_k^{\Omega^2 g}  (\Omega f) \stackrel g=\bar \Omega^{1-k} \delta_k^g f$.
\\
\item\label{OS} {\it On-shell property:} $\delta^g_1 \csigma-1=0=\delta_k \csigma$, when $k\in \{2,\ldots, d\}$ and $g_+^{\rm SY}=\csigma^{-2} \cg$ is the  singular Yamabe metric.
\end{enumerate}
\end{problem}
\noindent
In dimension $d$, the operator $\delta_{d+1}$ is a critical case because this is the order 
where one encounters the Neumann data for the singular Yamabe problem; see Equation~\nn{C0} and Section~\ref{DN}.

A  general theory of boundary operators
for general weights has been developed in earlier work of Gover and Peterson~\cite{GPt}. That work
mostly solves Problem~\ref{prob}.
Here we will give the first example of an explicit  fifth order operator. 
We  also show how to suitably restrict the underlying geometry in cases where the Gover--Peterson operators cannot be  applied. This relies on a subtle analysis of dimension-dependent poles in the operators of~\cite{GPt}.
In the case when  operators solving Problem~\ref{prob} exist, we  can construct scalar analogs of conformal fundamental forms.
  In particular, if $g_+ = \ceta^{-2} \cg$ is a CC metric, then 
$$
\pi_k(g_+) := \delta_k \ceta
$$
defines a $k^{\rm th}$ order {\it CC boundary curvature scalar}. 

\smallskip
By construction $\pi_1(g_+^{\rm SY})=1$ and 
\begin{equation}\label{pik0}
\pi_k(g_+^{\rm SY})=0
\end{equation}
for  all $k\geq 2$ such that the operator $\delta_k$ exists. As a consequence one can employ  {\it CC boundary curvature scalars} to compute expansion
coefficients for  asymptotic solutions to the singular Yamabe problem.
Details are given in  Section~\ref{EC}, while explicit formul\ae\ are given in Section~\ref{ER}.

The volume of a conformally compact manifold $(M_+,g_+)$ is ill-defined thanks to the metric's singular behavior along the boundary. However, the asymptotic behaviors of the volumes of regions when a collar neighborhood of the boundary is removed can be fruitfully studied; see in  particular~\cite{HS,GrVol}. Of particular interest is whether a finite, invariant ``renormalized volume'' can be extracted from these volume asymptotics.
Graham and Gursky have studied renormalized volumes of $(M_+,g_+^{\rm SY})$ with respect to the singular Yamabe metric~\cite{GG} for four-manifolds. In particular they find that when the conformal embedding $\Sigma\hookrightarrow(M,\cg)$ is umbilic, there exists a renormalized-volume type invariant. Renormalized volumes are typically  important global invariants associated to conformal manifolds with boundary. An interesting, but open problem, is whether singular Yamabe Dirichlet-to-Neumann maps  can capture the functional gradient of renormalized volume with respect to  variations of the boundary embedding.

For physical applications to general relativity, one is primarily interested in the case of 
spacelike boundaries~$\Sigma$ and positive scalar curvature. Indeed, upon adjusting for the change of metric signature, our results can  be applied to the study of asymptotically de Sitter spacetimes. In particular they can be used to complete the characterization of conformal infinities of asymptotically de Sitter spacetimes
given in~\cite{GoKo}. Another more or less direct application is to the characterization of solutions to  Kazdan--Warner-type problems~\cite{KW} generalizing the singular Yamabe problem to the case of prescribed (but not necessarily constant) scalar curvature.

%
%
%
%
%
%
%
%
%

%
%
%

\subsubsection*{Remembering Stanley Deser}

\begin{center}  By Andrew Waldron\end{center}

\noindent
Stanley Deser was a truly wonderful human being. He had a highly developed sense of humor---everybody walked away smiling after spending any amount of time with him. He was also highly knowledgeable and incredibly intelligent. Chances are you also walked away wiser. I first met Stanley via a series of terse emails leading to a, for me formative,  post-doctoral position at Brandeis in 1999. After leaving Brandeis, Stanley and I continued our highly fruitful collaboration over Skype. 
Those meetings followed Stanley's standard---50\% jocular persiflage plus 50\% question bombardment---format.
Back at Brandeis, Stanley usually arrived at work shortly after lunch, oftentimes muttering about seemingly agrarian peculiarities such as Bach and Cotton tensors. Later I realized this was because he was deeply interested in conformal anomalies. Indeed his work with Adam Schimmer on such was not only fundamental for the modern day understanding of how physical theories behave at different scales, but also laid the foundation for many mathematical studies of conformal manifolds. 
Curiously enough, my work with Stanley on partially massless higher spin theories, led me first to Tom Branson and in turn Rod Gover, another long term collaborator and expert in all things conformal. In that sense there is direct line from those early mutterings to the current paper; indeed the operator in Equation~\nn{PM} is the conformal antecedent of the spin~2 partially massless theory. 
I am sure Stanley would have enjoyed the results presented here and sorely miss being able to discuss them with him. Maybe a discussion would have gone something  like this (from one of our last email exchanges):

\begin{quote}
{\it A:} ${\tt What\ about\ a\ prequel\  to\  the\ Bible?\ The\
 last\ line\ reads}$
 
  ${\tt \qquad\qquad\quad ``in\ the\ end\ there\ was\ light".}$\\[-3mm]

{\it S:} ${\tt that's\ a\ typo:\ in\ the\ end\ there\ was\ gravity!}$
\end{quote}

\subsection{Conventions}\label{conv}
We denote by $\nabla^g_a$ the Levi--Civita connection of a metric $g$ in an abstract index notation (see~\cite{Penrose}), while $\nabla_v$ denotes covariant derivatives in the direction of any vector $v$---the superscript $g$  is dropped whenever context allows.
Our curvature convention is 
$$[\nabla_a,\nabla_b]v^c =:R_{ab}{}^c{}_d v^d$$ while the Ricci and scalar curvature tensors are respectively given by
$$Ric_{ab}:=R_{ca}{}^c{}_b\, ,\qquad
Sc:=g^{ab}R_{ab}\, .$$
The Schouten tensor $P_{ab}$  is defined, when $d\geq 3$,  by 
$$Ric =: (d-2) P + g J\, ,$$ where\ $J$ is the trace of $P$. 
Also note that when $d=2$, we define $J:=Sc/2$.

We will use a $\circ$ notation to denote a trace-free condition, for example $$\mathring X_{ab}:= X_{ab} -\tfrac 1d g_{ab} g^{cd} X_{cd}\in \Gamma(\odot_\circ^2T^*M)\, .$$
Clearly projecting to the trace-free part of a tensor does not depend on the choice of conformally related metrics $g\in \cg$.

Given a metric dependent tensor $T^g$, we define its weight $w$ linearized conformal variation by
  $${\updelta}^{\rm L}_w T:=\frac{\ext T^{e^{2t\varpi}g}}{\ext t} \big|_{t=0}-w\varpi T^g\, ,$$
  where $\log \Omega:=\varpi \in C^\infty M$.
An important result of~\cite{linear} is that is suffices to check that $\updelta^{\rm L}_w T=0$ in order to verify $T^{\Omega^2 g}=\Omega^w T^{g}$. 
Similarly, for tensors $\bar T^g$ defined along a hypersurface
we define
$$\bar {\updelta}^{\rm L}_w \bar T:=\frac{\ext \bar T^{e^{2t\varpi}g}}{\ext t} \big|_{t=0}-w\bar \varpi \bar T^g\, ,$$
where $\bar \varpi :=\varpi|_\Sigma$.

\section{Conformal Fundamental Forms}
\label{CFFs}

Given a conformal $d$-manifold $(M,\cg)$ with dimension $\bar d:=d-1$ boundary $\Sigma=\partial M$ 
one wishes to study the 
 {\it conformal embedding}
$$\Sigma \hookrightarrow (M,\cg)\, .$$
For that a key tool 
is  a sequence of conformally invariant
generalizations of the extrinsic curvature/second fundamental form~\cite{BGW}.
Firstly, given  metrics $g, \Omega^2 g\in \cg$, observe that the corresponding unit conormals are related by
$$
\hat n^{\Omega^2 g} = \bar \Omega\hh \hat n^g\, ,
$$
so this defines a section $\hat {\boldsymbol{n}}\in  \Gamma(T^*M[1]|_\Sigma)$. This is our first example of a {\it conformal hypersurface invariant}, meaning a suitably natural (see~\cite{CAG} for details)  tensor along~$\Sigma$, defined from the data $\Sigma\hookrightarrow (M,\cg)$, that is invariant under conformal rescalings of the bulk metric~$g$.

The induced metric $\bar g_{ab}\stackrel\Sigma = g_{ab}|_\Sigma-\hat n_a \hat n_b$ determines a weight $2$ conformal density-valued rank two tensor $\bar \cg_{ab}\in \Gamma(\odot^2T^*\Sigma[2])$, that in this context may be termed the {\it first conformal fundamental form}. Of course, 
we may also view this as the conformal class of metrics $\bar \cg$ induced along $\Sigma$. Note that there is no need to introduce new indices for hypersurface quantities once one remembers that the image $(TM|_\Sigma)^\top$ of the projection of $TM|_\Sigma$ 
to directions tangent to $\Sigma$ 
is isomorphic to $T\Sigma$.

\smallskip

The second conformal fundamental form~$\cIIo_{ab} \in \Gamma(\odot^2_\circ T^*\Sigma[1])$ has weight $1$ and is given by the trace-free part of the extrinsic curvature $\II_{ab}$,
$$
\IIo_{ab}\stackrel{g}{:=}\II_{ab}-\tfrac1{d-1}\bar g_{ab} \bar g^{cd} \II_{cd}\, .
$$
Recall that the extrisinc curvature $\II_{ab}:=\nabla^\top_a \hat n^{\rm ext}_b|_\Sigma$, 
where $\hat n^{\rm ext}$ is any smooth extension of the unit normal to $M$ and note that $\top$ is used to denote projection to  directions tangent to $\Sigma$. Also $H:=\tfrac1{d-1} \II_a{}^a$ is the mean curvature of $\Sigma$ with respect to the metric $g$ determining $\nabla$. 
 The scalar  $K:=\IIo^{ab} \IIo_{ab}$ will be termed the {\it rigidity density}, partly because it was employed by Polyakov to describe rigid string world sheets~\cite{Poly1}. 
Embeddings such that the second conformal fundamental form vanishes are termed {\it umbilic}. 

A third conformal fundamental form $\cIIIo_{ab}\in \Gamma(\odot_\circ^2T^*\Sigma[0])$ is given, in dimensions~$d>3$, by 
$$
\IIIo_{ab}
\stackrel g{:=}W_{a\hat n b\hat n}+\IIo^2_{(ab)_\circ}=(d-3)\big(P^\top_{(ab)_\circ}-\bar P_{(ab)_\circ} +H\IIo_{ab}\big)\, .
$$
The above-displayed equality follows from the trace of the Gauss equation. When $d=3$, we have that $W_{a\hat n b\hat n}+\IIo^2_{(ab)_\circ}\equiv 0$, so for our current purposes we  define  $\IIIo_{ab}:=0$ in that case. 
The third fundamental form gives a conformally invariant measure of the difference between bulk and boundary Ricci tensors.
Some authors~\cite{Vyatkin} employ to this end  the closely related conformally invariant {\it Fialkow  tensor}~\cite{Fialkow} $$F_{ab}:=P^\top_{ab}-\bar P_{ab} +H\IIo_{ab}+\tfrac12 H^2 \bar g_{ab}\, .$$  We shall term  $d\geq 4$ embeddings where $\boldsymbol F_{ab}=0$,  {\it Fialkow flat}.

Any fourth fundamental form 
$\cIVo_{ab}\in \Gamma(\odot_\circ^2T^*\Sigma[-1])$
has weight $-1$. A fourth fundamental form is  given by
$$
\IVo_{ab}\stackrel{g}{:=}
C_{\hat  n(ab)}^\top  -H W_{a\hat n b \hat n}+\tfrac1{d-5}\,  \bar \nabla^c W_{c(ab) \hat n}^\top=:
\widehat{\IVo\hh\hh}{\!}_{ab}+\tfrac1{d-5}\,  \bar \nabla^c W_{c(ab) \hat n}^\top\, .
$$
Note that we sometimes employ a hat notation for versions of tensors/operators where dimension-dependent poles have been removed.

In dimension four, if the conformal class $\cg$ admits a Poincar\'e--Einstein metric, it follows that both $\IIo$ and $\IIIo$ vanish, see~\cite{LeBrun, Goal} and~\cite{BGW}.
In this sense, these fundamental forms are obstructions to solving for an Einstein metric. In four dimensions, the fourth fundamental form plays a distinguished {\it r\^ole}: for Poincar\'e--Einstein metrics it computes both the image of the Dirichlet-to-Neumann map for the Einstein problem~\cite{BGKW} and the variational gradient of the renormalized volume of the Poincar\'e--Einstein manifold~\cite{Anderson}. More recently, also in four dimensions, fundamental forms were used to characterize the allowed stress energy of  asymptotically de Sitter spacetimes as well as initial, isotropic, cosmological singularities~\cite{GoKo,GoKoWa}.

\section{Conformally Compact Boundary Curvatures} \label{CCBC}

To construct CC boundary curvatures we first need to solve Problem~\ref{prob}. For that we recall the Laplace--Robin operator 
(see~\cite{GW}) specialized to 
 act on conformal densities.
Given any defining density $\ceta$, the  conformally invariant Laplace--Robin operator mapping
$${\mathscr L_{\ceta}}: \Gamma({\mathcal E} M[w])\to \Gamma({\mathcal E} M[w-1])$$
is defined for a choice of $g\in \cg$ and $\boldsymbol{\varphi}\in   \Gamma({\mathcal E} M[w])$ by
\begin{equation}\label{G&T}
\mathscr L_{\ceta} {\varphi}\stackrel{g}{:=}
(d+2w-2)\big(\nabla_\nu + w \rho_\eta) \varphi- \eta (\Delta^g+w J^g)\varphi
\, .
\end{equation}
In  the above, the conormal $\nu:=\ext \eta$ and $\rho_\eta:=-\tfrac1d (\Delta^g \eta + J^g\eta)$.
Specializing to the case where $\boldsymbol\varphi$ is the weight $w=1$ defining density $\ceta$ itself, we can recover the scalar curvature of the CC metric
$g_+=\ceta^{-2} \cg$, via
\begin{equation}\label{Sctr}
{\mathscr L}_{\ceta} \ceta = -\hh\frac{Sc^{g_+}}{d-1}\, .
\end{equation}
The above can be easily verified by examining the conformal transformation rule for scalar curvatures. Moreover
\begin{equation}\label{I2}
\frac{1}d\hh{\mathscr L}_{\ceta} \ceta\stackrel{g}=
|\nu|^2_g +2 \rho_\eta \eta
\, ,
\end{equation}
measures the squared-length of the conormal $\nu$ along $\Sigma$, where $\eta=0$. Also, along~$\Sigma$, the Laplace--Robin operator becomes a conformally invariant Robin-type normal operator with leading term $\nabla_\nu$. Indeed, it can be  checked (see~\cite{GPt}) that, along~$\Sigma$,
$$
\mathscr L_{\ceta}^k \, \pdot \, \big|_\Sigma = \prod_{i=1}^k (d+2w-k-i) \nabla^k_\nu\, \pdot \,\big|_\Sigma  + {\rm ltots}\, ,
$$
where once again ``ltots'' denotes terms of lower transverse derivative order, meaning lower order in derivatives normal to $\Sigma$.
Evidently from Equation~\nn{Sctr}, specializing to a singular Yamabe metric $g_+^{\rm SY} = \csigma^{-2}\cg$, then $\mathscr L_{\csigma}\csigma=1$ and for any $k\geq 2$ one has 
$$
\mathscr L_{\csigma}^k \csigma=0\, .
$$
Hence, Problem~\ref{prob} is solved by setting
$$
\delta_k=\frac{\mathscr L_{\ceta}^k\,\pdot \hh \big|_\Sigma}{ \prod_{i=1}^k (d+2w-k-i) }\, ,
$$
so long as the denominator is non-zero. 
Indeed, the precise structure of dimension-dependent poles in such operators is important.
Because we are here acting on weight one densities the denominator is
$$
(d-k+1)(d-k)\cdots(d-2k+2)\, .
$$
Observe that the Robin-type operator mapping, for some choice of $g\in\cg$,
\begin{equation}\label{Littleredbirds}
\boldsymbol{\varphi}\to (\nabla_\nu +w \rho_\nu)\varphi\big|_\Sigma\, ,
\end{equation}
is conformally invariant, so from Equation~\nn{G&T} we see that $(d-2k+2)^{-1} \mathscr L_{\ceta}^k\hh\pdot\hh \big|$ is both well-defined and conformally invariant when $d=2(k-1)$. In fact, there are other removable singularities when the operator $\mathscr L_{\ceta}^k$ is suitably improved by adding conformally invariant subleading transverse order terms, see Gover and Peterson~\cite{GPt}. The pole at $d=k-1$, cannot be removed this way. While such modified operators can be used to solve Problem~\ref{prob} for further values of $k$, they cannot be used when $k=d+1$. Indeed, from Equation~\nn{C0}, we see that this is precisely the order where, in general,  the singular Yamabe solution $\csigma$ fails to be smooth along $\Sigma$ so that $\nabla_\nu^{d+1} \sigma$ is not well-defined along $\Sigma$.
 In  $d$ dimensions, for a choice of $g\in \cg$ we will explicitly construct operators (for cases $k\leq 5$ and motivated by the methods of~\cite{GPt}), such that
\begin{equation}\label{ren}
\delta_{k} \propto \frac{\bar B_{k-1}}{d-k+1}+\delta_k^{\rm ren}\, ,
\end{equation}
where the operator $\delta_k^{\rm ren}$ is not singular when evaluated in $d=k-1$ dimensions and~$\bar B_{k-1}$ is some tensor defined in any dimension that specializes to  the conformally invariant obstruction density for the 
singular Yamabe problem in 
$d=k-1$ dimensions.  Upon suitably specializing the embedding geometry, the operator $\delta_{d+1}^{\rm ren}$ also yields an invariant that characterizes a Dirichlet-to-Neumann map for the singular Yamabe problem. 

In fact, there are further poles that cannot be removed by applying the Gover--Peterson improvement. There are  two  cases to consider.
The first is when  the residue of the pole is given by the restriction to the critical dimension $d_{\rm crit}$ of a tensor that is conformally invariant in all dimensions. Then one can  construct a conformally invariant operator by simply removing the invariant pole term. One then must check whether the onshell Property~(\ref{OS}) still holds.
In second case the residue is some tensor that has no invariant continuation to general dimensions. In that case the variation of the residue takes the form $(d-d_{\rm crit}) X$ for some tensor~$X$. This determines the failure of the remainder to be invariant after the pole is removed. Thus one must
impose the condition~$X=0$ 
in  order that the pole can be removed and still yield an invariant.
For example, when $d=3$ requiring $\pi_1=1$, or equivalently that the geometry is asymptotically hyperbolic,  suffices to remove such a residue this way.
 The fact that we can find CC boundary scalars whose invariance in certain critical dimensions is conditioned on restrictions upon  lower order CC boundary scalars 
mimics a similar feature observed in~\cite{BGW} for conformal fundamental forms.

\subsection{Expansion Coefficients}\label{EC}

To construct an asymptotic solution to the singular Yamabe problem, 
typically one considers a distinguished compactified metric $g_{\rm GL}(\bar g)\in \cg$, where $\bar g\in \cg $ denotes a given choice of  boundary  metric in the conformal class of metrics induced by $\cg$. For that  first observe that
given any CC metric $g_+=\ceta^{-2}\cg$, 
using Equations~\nn{Sctr} and~\nn{I2} we can canonically construct, at least in a collar neighborhood of $\Sigma$, an asymptotically hyperbolic metric $$g_+^{\rm AH}=\hat \ceta^{-2} \cg\, ,$$ 
by setting
$$
\hat \ceta =\sqrt{-\frac{d(d-1)}{Sc^{g_+}}}\: \ceta\, .
$$
Then one relies on a key result of Graham and Lee~\cite{GrahamLee} that given any asymptotically hyperbolic metric  $g_+^{\rm AH}=\hat\ceta^{-2} \cg$ as well as a boundary metric $\bar g\in \bar \cg$, there exists a ``Graham--Lee''  scale $\ctau_{\rm GL}$ such that, also at least in a sufficiently small collar neighborhood of $\Sigma$, the function 
$$
r:=\frac{\hat \ceta}{\ctau_{\rm GL}}
$$
is the geodesic distance to the boundary as measured by the compactified metric $g_{\rm GL}=\ctau_{\rm GL}^{-2} \cg$. Away from the collar there may be no distinguished choice for $\ctau_{\rm GL}$, but this ambiguity is irrelevant for our considerations. Now, given the Graham-Lee compactified metric $g_{\rm GL}(\bar g)$ for some (reference) CC metric $g_+$ and corresponding distance function $r$, one expands the singular Yamabe scale $\csigma$ expressed in the Graham--Lee choice of scale as
$$
\sigma\stackrel{g_{\rm GL}\!(\bar g)}= r(s_0 +  s_1 r +  s_2r^2 +\cdots)\, ,
$$
to as high order as regularity of the singular Yamabe scale $\csigma$ allows.
Here, one can require that the functions~$s_i$ obey the Lie derivative condition ${\mathcal L}_{\frac\partial{\partial r}}s_i=0$, since the Graham--Lee metric $g_{\rm GL}$ admits a normal form 
$$g_{\rm GL}(\bar g)= dr^2 + h\, ,$$ where the symmetric tensor $h(r)$ obeys
$h({\frac\partial{\partial r}},\hh\pdot\hh)=0$, so that $g^{-1}_{\rm GL}(dr,\hh \pdot\hh) = {\frac\partial{\partial r}}$. Moreover, since the unit conormal $\hat n$ with respect to $g_{\rm GL}$  equals $\ext r|_\Sigma$, it follows that acting on scalars, the operator 
$$
\hat n^{a_1}\cdots  \hat n^{a_k} \nabla_{a_1}\cdots \nabla_{a_k}\hh\pdot\,|_\Sigma={\mathcal L}_{\frac\partial{\partial r}}^k\hh\pdot\,|_\Sigma\, .
$$
Hence  from Equation~\nn{pik0} we have that $-1/k!$ times the terms labeled ``ltots'' in Point~(\ref{NO}) of Problem~\ref{prob} yields the expansion coefficient $s_{k-1}$
in the above expansion. Explicit results for $s_0,\ldots s_4$ are given in  Section~\nn{EC5}.

\section{The Dirichlet-to-Neumann Map}\label{DN}\label{bigbluebirds}

We wish to characterize the Neumann data for the singular Yamabe problem in a canonical way.
Given a second
 order elliptic PDE one may prescribe only its Dirichlet data to study its solution space. Thus, an apparent freedom to prescribe further Neumann data is only a property of  formal asymptotic solutions, for example the density $\boldsymbol{C}$ in Equation~\nn{C0} depends on the original choice of defining density~$\csigma_0$ so is not an absolute invariant of the underlying geometric structure. However, there are situations, such as Poincar\'e--Einstein metrics in even dimensions, where 
one can assign a conformally invariant Dirichlet-to-Neumann map to the geometric structure~\cite{LittleRobin,BGKW}. Rather than relying on  dimension parity-dependent properties of Poincar\'e--Einstein metrics, here we will construct 
a conformally invariant Dirichlet-to-Neumann map in the case when the extrinsic geometry of the conformal embedding $\Sigma\hookrightarrow(M,\cg)$ is suitably restricted. While the procedure is general in nature, we focus on cases where $d\leq 4$.

Our plan is to ``renormalize'' the operator $\delta_{k}$ in Equation~\nn{ren}, by replacing it with~$\delta_k^{\rm ren}$ in the case where $d=k-1$. For this we need to know whether $\delta_k^{\rm ren}$ is invariant. For that we compute the conformal variation of $\bar B_{k-1}$ in $d$ dimensions. 
In the current setting involving natural tensors, it can be shown (see for example~\cite{GPt}) that this must take the form
$$
\bar B^{\Omega^2 g}_{k-1}-\bar \Omega^{1-k}
\bar B^g_{k-1} = (d-k+1) \beta \, ,
$$
for some tensor $\beta$ that is non-singular in $d=k-1$ dimensions. 
We then analyze sufficient conditions on the extrinsic geometry of the conformal embedding $\Sigma\hookrightarrow(M,\cg)$ such that $\beta$ vanishes. When such conditions exist it follows, upon imposing them,  
 that the renormalized operator $\delta_{d+1}^{\rm ren}$ is conformally invariant. Dimensional continuation arguments of this type 
can be rather subtle, so we also
 check this invariance explicitly. 
 
 The operator $\delta_{d+1}^{\rm ren}$ can now be used to invariantly extract singular Yamabe Neumann data. However, since its leading transverse order is $d+1$, it is ill-defined acting on solutions with only $C^d$ regularity. Hence---see Equation~\nn{C0}---we must further restrict the extrinsic geometry such that the obstruction $\bar {\boldsymbol B}_d=0$.
For our current purposes let us call a conformal embedding such that both $\beta$ and $\bar {\boldsymbol B}_d$ vanish {\it superumbilic}, because in dimension $d=3$ this amounts to the embedding being umbilic. For superumbilic geometries, we can define a {\it Dirichlet-to-Neumann map} which takes as input a closed, oriented, superumbilic, hypersurface embedding in a suitable conformal manifold, 
and outputs a conformal density according to
$$
\big(\Sigma\hookrightarrow 
(M,\cg)\big) \longmapsto \delta^{\rm ren}_{d+1} \csigma_{\rm SY}=:\pi_{d+1}^{\rm ren}\in \Gamma({\mathcal E}\Sigma[-d])\, .
$$
In the above  $\csigma_{\rm SY}$ is the uniquely determined solution to the singular Yamabe problem defined on the interior region $M_+$ determined by the closed hypersurface $\Sigma$. We also require that the interior region $M_+$ defined by $\Sigma$ has itself boundary $\Sigma$.

At this point it might be unclear how the above map is related to an undetermined term in an asymptotic solution such as the density $\boldsymbol{C}$ in Equation~\nn{C0}. One may even think, based on Equation~\nn{pik0} and the discussion of Section~\ref{EC},   that the image $\pi_{d+1}^{\rm ren}$ of the Dirichlet-to-Neuman map vanishes.
Indeed, when the geometry is superumbilic,
the singular Yamabe scale is smooth and the quantity
$
{\mathscr L}_{\csigma}^{d+1} \csigma
$
is well-defined and vanishing. However, since $\pi_{d+1}^{\rm ren}$ is constructed via the dimensional continuation  coming from Equation~\nn{ren}, it does not {\it a priori} follow that $\pi_{d+1}^{\rm ren}=0$. In fact,  given  some $\bar g\in \bar \cg$,
in the Graham--Lee scale 
$\ctau_{\rm GL}^{\rm SY}$ 
of the singular Yamabe metric itself, for the case $d=3$,
we will find
$$
\pi^{\rm ren}_{4} 
\stackrel{g^{\sss\rm SY}_{\sss \rm GL}}= -3 \hat n^a \hat n^b \hat n^c (\nabla_{c} P_{ab}) \big|_\Sigma\, .
$$

When one is equipped with a reference CC metric $g_+=\ceta^{-2}\cg$   with geodesic distance function $\eta\stackrel{g_{\sss \rm GL}}=r$,  and the embedding is superumbilic, one has the smooth expansion
$$
\sigma\stackrel{g_{\sss \rm GL}}= r(s_0+ s_1 r + \cdots + s_{d-1} r^{d-1} + C r^d) + {\mathcal O}(r^{d+2})\, .
$$
Thus
$$
\pi_{d+1}^{\rm ren} = (d+1)! \, C + {\rm ltots}\, ,
$$
and so can be used to compute the formally undetermined coefficient $C$ in the asymptotic solution. Explicit formul\ae\ are given in Equations~(\ref{pi4ren}, \ref{pi5ren}) below.

\section{Explicit Results}\label{ER}

We are now ready to present  the first five CC boundary curvatures. The fifth is novel. As discussed in Section~\ref{CCBC} these can be computed by studying powers of the Laplace--Robin operator. In principle it would also be possible to employ the tractor calculus construction of~\cite{GPt}, but at fifth order this involves a rank five tractor, a difficulty avoided by the Laplace--Robin method. Another method is simply to write down all possible extrinsically coupled Riemannian operators of the correct homogeneity and then determine which linear combination of terms solves Problem~\ref{prob}.
For this it is sufficient to study linearized conformal variations; see~\cite{linear} and the discussion in Section~\ref{conv}. For the results that follow, we employed both the Laplace--Robin operator construction and linearized variation methods using the symbolic manipulation software FORM~\cite{Jos}. 
Note that, as mentioned earlier, one can also act with the normal operators of~\cite{GPt} on the scalar curvature of the CC metric $g_+=\ceta^{-2}\cg$. That method produces curvatures quadratic in $\ceta$ and returns zero when  $\ceta$ is the singular Yamabe scale.

\subsection{First conformally compact boundary scalar}

The operator $$\delta_1:\Gamma({\mathcal E}_0M[1])\to \Gamma({\mathcal E}M[0])\, ,$$
defined for some $g\in \cg$
by
$$
\delta_1  \stackrel{g}{:=} \hat n^a \nabla_a \hh \pdot \, |_\Sigma\, .
$$
Observe that acting on the space $\Gamma({\mathcal E}_0M[1])$, 
the  {Robin} operator of Equation~\nn{Littleredbirds} becomes the above conformally invariant Neumann-type operator.
 Hence~$\delta_1$ obeys Properties~(\ref{NO}) and~\nn{CW}. Now let $g_+=\ceta^{-2}\cg$ be CC. Then the first CC scalar curvature is simply
$$
\pi_1\stackrel{g}{:=}\delta_1 \eta = \hat n.\nu|_\Sigma\, .
$$
Here $\nu:=\ext \eta$.   
By virtue of Equation~\nn{I2}, when $\ceta$ equals the singular Yamabe defining density $\csigma$,
 it follows that $|n|_g^2\stackrel\Sigma=1$ so that $\hat n=n|_\Sigma$, where $n:=\ext \sigma$. Hence, in that case $\delta_1\sigma =1$ and the on-shell Property~\nn{OS} holds.
 
Notice that the condition $\pi_1=1$, is equivalent to requiring that $g_+$ be asymptotically hyperbolic. Moreover, if $\pi_1^{\rm ext}$ is any smooth extension of $\pi_1$ to $M$, then the metric
$$
g_+^{\rm AH} = \frac{g_+}{(\pi_1^{\rm ext})^2} 
$$ 
is necessarily asymptotically hyperbolic.

\subsection{Second conformally compact boundary scalar}
The operator 
$$\delta_2:\Gamma({\mathcal E}_0M[1])\to \Gamma({\mathcal E}M[-1])\, ,$$ is defined for some $g\in \cg$
by
\begin{equation}\label{delta2}
\delta_2 \stackrel g{ :=} \hat n^a \hat n^b\nabla_a\nabla_b \, \pdot\, |_\Sigma-H \delta_1  \, .
\end{equation}
Clearly property~\nn{NO} holds, 
and using the conformal transformation rule for mean curvatures
$$
\bar \Omega H^{\Omega^2g}= H^g + \hat n.\Upsilon|_\Sigma\, ,
$$
where $\Upsilon:=\ext \log \Omega$,
it is not difficult to verify property~\nn{CW}. For the onshell property~\nn{OS}, we need a short computation. Firstly,
it can easily be proved (see for example~\cite{Goal}) that 
$$
\rho_\sigma|_\Sigma \stackrel g= -H^g\, ,
$$
when $g_+^{\rm SY}=\csigma^{-2}\cg$ is a singular Yamabe metric. 
But, in that case (calling $\rho:=\rho_\sigma$) 
$$
\hat n^a \hat n^b\nabla_a\nabla_b \sigma|_\Sigma =
\hat n^a  n^b \nabla_a n_b|_\Sigma
=\frac12 \hat n^a \nabla_a n^2|_\Sigma = 
\frac12 \hat n^a \nabla_a (1-2\rho\sigma)|_\Sigma =-\rho|_\Sigma\, .
$$
In the above we repeatedly used Equation~\nn{I2}. The onshell property~\nn{OS} now follows. Hence, relaxing to the case when $g_+=\ceta^{-2}\cg$ is CC, 
$$
\pi_2 := \hat n^a \hat n^b \nabla_a \nu_b|_\Sigma- H \pi_1
$$
is the second CC boundary scalar.

\begin{remark}
The second order normal operator $\delta_2$ becomes critical in dimension  $d=1$. To see that, call $\kappa:=\bar g^{ab} \II_{ab}$. Then the mean curvature $H=\tfrac1{d-1} \kappa$, and we can rewrite Equation~\nn{delta2} as
$$
\delta_2   =\hat n^a \hat n^b\nabla_a\nabla_b \,\pdot\,|_\Sigma-\tfrac1{d-1} \kappa\,  \delta_1 \, .
$$
Although $\kappa=0$ in dimension $d=1$, it does not follow that $\hat \delta_2 :=\hat n^a \hat n^b\nabla_a\nabla_b \hh\pdot\,|_\Sigma$ is an invariant operator when $d=1$. Rather
$$
\bar \Omega\hh\hat  \delta_2^{\Omega^2 g} \circ\Omega  =\hat \delta_2  + \hat n^a \Upsilon_a\hh\pdot\,|_\Sigma\, .
$$
This is an example of the general phenomenon discussed in Section~\ref{bigbluebirds}.
\hfill$\blacklozenge$
\end{remark}

\subsection{Third conformally compact boundary scalar}\label{third}

At normal orders three and higher, another subtle phenomenon arises.
 Namely, conformally invariant operators and curvatures intrinsic to the boundary appear as residues of poles when the dimension~$d$ is viewed as a parameter. 
We will use the notation $\delta_k^{[I]}$ where $I$ is some set of positive integers  for which the operator is not defined for reasons of  poles $1/(d-x)$ for $x\in I$. 
Recall that from Equation~\nn{C0}, the order $k$ CC boundary scalar is sensitive to Neumann data in dimension $d=k-1$.
This will correspond to a pole $1/(d-k+1)$
in the operator $\delta_k^{[I]}$. We are less interested in poles $1/(d-x)$ for integers $x<k-1$, because these probe data beyond the order of the Neumann term. Thus we primarily focus on  integers $k-1\leq x\in I$. 

\medskip

The operator 
$$\delta_3^{[2,3]}:\Gamma({\mathcal E}_0M[1])\to \Gamma({\mathcal E}M[-2])\, ,$$ is defined for $d\neq 3$ and some $g\in \cg$
by
$$
\delta_3^{[2,3]}  \stackrel g{:=}\hat \delta_3 -\Big(\tfrac{3}{d-3} \bar \Delta 
-\tfrac{2}{d-2} \hh K\Big)\delta_1\,
\, ,$$
where 
$$
\hat \delta_3 \stackrel{g}{:=} \hat n^a  \hat n^b \hat n^c \nabla_a \nabla_b \nabla_c\,\pdot\,|_\Sigma
+2 P_{\hat n\hat n} \hh\delta_1
\,\,  .
$$
For later use, in dimensions $d>2$ and for any $x\neq 2,3$, we also define the operator
$$
\delta_3^{(x)}
\stackrel{g}{:=}
\hat n^a  \hat n^b \hat n^c \nabla_a \nabla_b \nabla_c \hh\pdot\,|_\Sigma
+\Big(2 P_{\hat n\hat n} 
+\tfrac{2}{x-2} \hh K\Big)\delta_1
 -\tfrac{3}{x-3} \bar \Delta 
\delta_1 
\, .
$$

Observe that the singularity of $(d-3)\delta^{[3]}_3$ at $d=3$ is removable and yields the  operator $-3\bar\Delta \delta_1$ which is conformally invariant precisely in dimension $d=3$ because the Laplacian $\bar \Delta$ intrinsic to~$\Sigma$  is conformally invariant acting on functions, or equivalently conformal densities of weight $0$. (There is also a removable singularity in $(d-2)\delta^{[3]}_3$ when $d=2$, but the corresponding conformal invariant $K:=\IIo^{ab}\hh\IIo_{ab} $ vanishes because the second fundamental form is pure trace for curves embedded in two-manifolds.)
That the operator $\delta^{[3]}_3$ obeys Properties~(\ref{NO}-\ref{OS}) in general dimensions~$d$ can be checked relatively easily by direct computation. 

\medskip
In  dimensions $d\geq 4$, we now obtain a third 
CC boundary scalar
\begin{equation}\label{pi34}
\pi_3\stackrel{g}{:=} 
\hat n^a  \hat n^b \hat n^c \nabla_a \nabla_b \nu_c|_\Sigma+\big(2 P_{\hat n\hat n} 
-\tfrac{3}{d-3} \bar \Delta 
+\color{black}\tfrac{2}{d-2} K \big)\pi_1
\, .\end{equation}

Now let us consider the case $d=3$.
The linearized conformal variation of $\bar \Delta \pi_1$ is
$$
\bar \updelta^L \bar \Delta \pi =(d-3) \nabla_{\bar \Upsilon} \pi_1\, ,
$$
where $\bar \Upsilon := \ext \bar \varpi$.
Thus, as explained in Section~\ref{bigbluebirds}, we remove the term $-\tfrac{3}{d-3} \bar \Delta$ from~$\pi_3$ and then continue to $d=3$ dimensions and require that
 $\pi_1$ is constant, in which case also
$\bar \Delta \pi_1=0$. 
Then, as can be easily verified,
$$
\pi_3|_{
\pi_1=1,d=3
}:= 
\hat n^a  \hat n^b \hat n^c \nabla_a \nabla_b \nu_c|_\Sigma+2 (P_{\hat n\hat n} +K)
\, ,$$
is conformally invariant meaning $\pi_3|_{
\pi_1=1,d=3
}(\Omega^2 g,\Omega\sigma) = \bar\Omega^{-2}\pi_3|_{
\pi_1=1,d=3
}(g,\sigma)$.
Recall that for  asymptotically hyperbolic metrics $\pi_1=1$.
So $\pi_3|_{
\pi_1=1,d=3
}$ is invariant  when $g_+^{\rm AH}= \ceta^{-2}g$ is asymptotically hyperbolic.
This is the conditional invariance phenomenon first discussed in Section~\ref{CCBC}.

There is a rather  interesting way to relax the 
asymptotically hyperbolic condition for embeddings that are nowhere umbilic. For that we define an operator mapping $\Gamma({\mathcal E}_0M[1])\to \Gamma({\mathcal E}[-4])$ for $d>2$ and any $g\in \cg$ by
\begin{equation}\label{Kbowtie}
K\bowtie \hat \delta_3  \stackrel g{:=} K \big(\hat n^a \hat n^b \hat n^c \nabla_a \nabla_b \nabla_c  \, \pdot\,|_\Sigma+2P_{\hat n\hat n}\hat n^a\nabla_a \, \pdot \, |_\Sigma\big) +\tfrac32 (\bar \nabla^a K) \bar \nabla_a \, \delta_1\, .
\end{equation}
Since the rigidity density vanishes iff the embedding is umbilic, the above operator 
can be used to define a third CC boundary scalar for $d=3$, nowhere-umbilic embeddings,
$$
\pi_3^{\scalebox{.7}{$\mbox{\sout{ umbilic }}$}}\stackrel g{:=}\hat n^a \hat n^b \hat n^c \nabla_a \nabla_b \nu_c |_\Sigma+\big(2P_{\hat n\hat n}+\tfrac32 (\bar \nabla^a \log K) \hh\bar \nabla_a \big) \pi_1\, .
$$

\smallskip

The situation in $d=2$ dimensions is particularly subtle. 
The pole term $\frac{2}{d-2}K$ in Equation~\nn{pi34} poses no difficulty
because it is invariant in general dimensions and moreover $\IIo$ and hence $K$ vanish identically in  two dimensions. The issue is instead the Schouten tensor term, since in $d$ dimensions
$P^g=\frac1{d-2}\big(Ric-\frac1{d-1}\hh g\hh\hh Sc \big)$, and there is no condition on the geometry alone for which the linearized conformal variation
$-(d-2)\nabla \Upsilon$ of the tensor $Ric-\frac1{d-1}g\hh Sc$ vanishes.
However an operator of transverse order three ought access the undetermined Neumann data for  singular Yamabe solutions when $d=2$. Therefore we now no longer search for an invariant third order CC boundary scalar that vanishes when evaluated on singular Yamabe metrics, but instead search for a  third normal order operator  with a distinguished conformal transformation. 
For that we observe that in $d$ dimensions one has the Ricci-type identity
$$
\hat n^a \hat n^b P_{ab}|_\Sigma = J|_\Sigma -\bar J +\tfrac{d-1}2H^2-\tfrac{1}{2(d-2)} K\, .
$$
Hence, in $d$ dimensions, we see from Equation~\nn{pi34} that the scalar $$(\hat n^a  \hat n^b \hat n^c \nabla_a \nabla_b \nu_c+2 J)|_\Sigma-2\bar J +(d-1)H^2$$ is conformally invariant. However, because there is no intrinsic scalar curvature quantity in one dimension, $\bar J$ is not defined when $d=2$. Hence let us consider
$$
\pi_3^{\rm rel}\stackrel g{:=}(\hat n^a  \hat n^b \hat n^c \nabla_a \nabla_b \nu_c+2 J)|_\Sigma+H^2\, .
$$
We note that this defines a {\it relative conformal invariant}, meaning that if ${\mathcal J}=[\bar g,j]$ is any equivalence class of boundary tensors such that$$
(\bar g,j)\sim \Big(\bar \Omega^2 \bar g, \bar \Omega^{-2}\big(j-\bar \Delta \log \bar \Omega+\tfrac12 |\ext \log \bar \Omega|^2_{\bar g}\big)\Big)\, ,
$$
then the combination
$$
\pi_3^{\rm rel}-2{\mathcal J} 
$$
is conformally invariant. The additional structure ${\mathcal J}$ required along $\Sigma$ amounts to a tensor that transforms in the same way as $\bar J$ would have had it existed. This can be viewed as a one dimensional analog of a M\"obius structure.

\subsection{Fourth conformally compact boundary scalar}

The fourth CC boundary scalar measures local deviations of CC metrics from the singular Yamabe metric  in dimensions $d>3$. In three dimensions, we can use it to extract the Willmore invariant as well as study Neumann data for the singular Yamabe metric. (As promised, from hereon we focus on operators in dimensions $d>2$ so will need no special treatment of~$d=2$ poles.) Before writing down the requisite normal operator~$\delta_4^{[3,5]}$, we need some additional technology.

Let $\mathring{\bar X}_{ab}$ be a symmetric trace-free tensor on $\Sigma$. Then we can define an operator~${\sf L}_{ab}$, mapping $\Gamma(\odot_\circ^2T^*\Sigma[4-d])\to 
\Gamma({\mathcal E}\Sigma[-d])$,
 such that
$$
\bar \Omega^{d}{\sf L}_{ab}^{\Omega^2 g}(\bar \Omega^{4-d}\mathring{\bar X}^{ab})={ {\sf L}}^g_{ab}\mathring{\bar X}^{ab}\, ,
$$
given for any $g\in \cg$ by
\begin{equation}\label{Lop}
{\sf L}^{ab} \mathring{\bar X}_{ab}\stackrel{g}{:=}
\left\{\begin{array}{cl}\bar \nabla^a \bar \nabla^b \mathring {\bar X}_{ab} + \bar P^{ab} \mathring{\bar X}_{ab}\, ,&d>3\, ,
\\[1mm]
\bar \nabla^a \bar \nabla^b \mathring {\bar X}_{ab} +  P^\top_{ab} \mathring{\bar X}^{ab}
+H \IIo^{ab} \mathring{\bar X}_{ab}\, ,&d=3\, .
\end{array}\right.
\end{equation}
Also, in any dimension $d$,  let us call ${\sf L}^{(3)}_{ab}\mathring{\bar X}^{ab}:=\bar \nabla^a \bar \nabla^b \mathring {\bar X}_{ab} +  P^\top_{ab} \mathring{\bar X}^{ab}
+H \IIo^{ab} \mathring{\bar X}_{ab}$.
Of course, this quantity only enjoys a conformal invariance property when $d=3$.
 In particular, let us define in any dimension the scalar
 $$
\bar B_{3}\stackrel{g}{ :=} -\tfrac13 {\sf L}^{(3)}_{ab} \IIo^{ab}\, ,
$$
which equals the Willmore invariant  when $d=3$.
Its ($d$-dimensional)  linearized conformal variation is given by 
\begin{equation}\label{slugmister}
\bar{\updelta}^L_{-3}\bar B_{3}=(d-3)\big(\IIo^{ab} \bar \nabla_a \bar \Upsilon_b + 2 (\bar \nabla^a \IIo_{ab})\bar \Upsilon^b\big)=
(d-3)\hh\bar  \square^{\IIo}_{(3)} \bar \varpi
\, ,\end{equation}
where  the operator
 $\bar \square^{\IIo}_{(d)}\stackrel g{:=}\big((d-2)\hh
 \IIo^{ab}\bar\nabla_a\bar\nabla_b
+
2(\bar\nabla_a\hh  \IIo^{ab})\bar \nabla_b
\big)$ 
is conformally invariant
acting on functions  in~$d$ dimensions.

In general dimensions we employ the notation $\bar {\sf L}^{ab} \mathring{\bar X}_{ab}\stackrel g{:=} \bar \nabla^a \bar \nabla^b \mathring {\bar X}_{ab} + \bar P^{ab} \mathring{\bar X}_{ab}$ even in  the case when the tensor $\mathring{\bar X}_{ab}$ may not have the correct transformation property/weight required for this combination to be invariant.

We need one further operator intrinsic to the conformal manifold $(\Sigma, \bar \cg)$:
 Let $\bar f$ be a smooth function on $\Sigma$. Then  when $\bar d= {\rm dim}\Sigma\geq 2$, we define the {\it Yamabe operator}~$\bar \square$, mapping $\Gamma({\mathcal E}\Sigma[1-\frac {\bar d}2])\to \Gamma({\mathcal E} \Sigma[-1-\frac{\bar d}2])$,
for any~$\bar g\in \bar \cg$ by
$$
\bar\square  \stackrel{\bar g}{:=} \bar \Delta   + \big(1-\tfrac {\bar d}2\big ) \bar J   \, .
$$ 
This obeys the conformal invariance property
$$
\bar \Omega^{1+\tfrac {\bar d}2}\,  \bar \square^{\bar \Omega^2 \bar g} \big(\bar \Omega^{1-\tfrac {\bar d}2}\bar f\big) = \bar \square^{\bar g}\bar f\, .
$$
Given an integer $x$,  let us denote in any dimension  the operator $\bar \square^{(x)}\stackrel{\bar g}{:=}\bar \Delta + \big(1-\tfrac{x}2\big)\bar J$. Note that, acting on a section $\bar f$ of ${\mathcal E}\Sigma[w]$, its linearized conformal variation is
\begin{equation}\label{boxshift}
\updelta^{\rm L}_{w} (\bar \square^{(x)}f) = (\bar d+2w -2) \bar \nabla_{\bar \Upsilon} f + (w-1+\tfrac x2) (\bar \nabla^a \bar \Upsilon_a) f\, . 
\end{equation}

\smallskip

We are now ready to build the fourth CC boundary curvature.
The operator 
$$\delta_4^{[3,5]}:\Gamma({\mathcal E}_0M[1])\to \Gamma({\mathcal E}M[-3])\, ,$$ is defined for $d\neq 3,5$ and some $g\in \cg$
by
\begin{equation*}
 \delta_4^{[3,5]}
 \stackrel{g}{:=}\hat \delta_4 
-\tfrac{6}{d-5}\bar \square^{(4)}\delta_2
-\tfrac3{d-3} \big(2\IIo^{ab}\IIIo_{ab}-3\bar B_{3}\big) \delta_1
\, ,
\end{equation*}
where
\begin{align*}
 \hat \delta_4
 &
 \stackrel{g}{:=}
 \hat n^a  \hat n^b \hat n^c \hat n^d \nabla_a 
\nabla_b \nabla_c \nabla_d\,\pdot\, \big|_\Sigma
+2 H \hat  \delta_3
+2 P_{\hat n\hat n} \delta_2
\\&\:
-\big(
8\hh  \IIo^{ab}\bar\nabla_a \bar \nabla_b
-\tfrac{3}{d-2} (\bar\nabla^a\bar \nabla^b \IIo_{ab})
+4 (\bar \nabla^a H)\bar \nabla_a
-2HP_{\hat n\hat n}
-3
\hat n^a \hat n^b \hat n^c (\nabla_a P_{bc })
\big)\delta_1
\, .
\end{align*}
Observe that the residues of the $d=3,5$ poles in $\delta_3^{[3,5]}$ are conformally invariant when $d=3,5$, respectively.

In dimensions $5\neq d\geq 4$ we now obtain a fourth CC boundary curvature
\begin{align*}
\pi_4&\stackrel{g}{:=} 
 \hat n^a  \hat n^b \hat n^c \hat n^d \nabla_a 
\nabla_b \nabla_c \nu_d|_\Sigma\!
-\tfrac{6}{d-5}\bar \square^{(4)}\pi_2\!
-\tfrac3{d-3} \big(2\IIo^{ab}\IIIo_{ab}\!-3\bar B_{3}\big) \pi_1
-\!\tfrac{4}{d-2} \!\hh H K \pi_1
\\&\:\:+\!2H \big(\pi_3\!+\!\tfrac3{d-3} \hh \bar \Delta \pi_1\big)
\!+\!2 P_{\hat n\hat n} \pi_2
\\&
\:\:-\Big(
8\hh  \IIo^{ab}\bar\nabla_a\! \bar \nabla_b
-\tfrac{3}{d-2} (\bar\nabla^a\bar \nabla^b \IIo_{ab})
\!\!+\!4 (\bar \nabla^a\! H)\bar \nabla_a
\!\!-\!2H\!\hh P_{\hat n\hat n}
\!\!-\!3\hat n^a \hat n^b (\nabla_{\hat n}\!\hh P_{ab })
\Big)\pi_1  .
\end{align*}
Let us now analyze this result. First note from Equation~\nn{pi34}, that the combination
$\pi_3\!+\!\tfrac3{d-3} \hh \bar \Delta \pi_1$
has no pole at $d=3$. The remaining singularity in $(d-3)\pi_4$ at $d=3$ is removable and, as promised,  yields the multiple 
$
9 \bar B_3 \hh \pi_1
$
of  the higher Willmore invariant~$\bar B_3$.

Next we consider the case $d=5$. Importantly, while $\bar \square^{(4)}\pi_2$ is invariant when $d=5$, its failure to be invariant in general dimensions is proportional to $d-5$; see Equation~\nn{boxshift}. Hence, along exactly the same lines as discussed in Section~\ref{third} for $\pi_3$ in $d=3$ dimensions, to obtain a fourth CC boundary curvature in $d=5$ dimensions, we must condition on vanishing of the second CC boundary curvature $\pi_2$, yielding the conformal invariant
\begin{align*}
\pi_4\big|_{\pi_2=0,d=5} &\stackrel{g}{:=}
 \hat n^a  \hat n^b \hat n^c \hat n^d \nabla_a 
\nabla_b \nabla_c \nu_d|_\Sigma +2H\pi_3
-\tfrac{15}2\IIo^{ab}\IIIo_{ab} \pi_1
\\&\: 
-\Big(
8\hh  \IIo^{ab}\bar\nabla_a\! \bar \nabla_b
+\tfrac12 
(\bar\nabla_a \bar\nabla_b \IIo^{ab}) + \tfrac32\bar P_{ab} \IIo^{ab}
+3H\bar \Delta
+\!4 (\bar \nabla^a\! H)\bar \nabla_a+\tfrac43 HK
\\&\qquad
-\, 2H\!\hh P_{\hat n\hat n}
-3\hat n^a \hat n^b (\nabla_{\hat n}\!\hh P_{ab })
\Big)\pi_1  .
\end{align*}

It now remains to study the case $d=3$. We already know that $(d-3) \pi_4$ recovers the higher Willmore invariant~$\bar B_3$, which obstructs smoothness of singular Yamabe solutions. We now wish to extract their Neumann data. We must consider the residue of the~$d=3$ pole in $\pi_4(g_+^{\rm SY})$ which is given by $
9 \bar B_{3}
$. Notice here that not only does~$\IIo^{ab}\IIIo_{ab}$ vanish when $d=3$ but more importantly, this tensor is invariant in general dimensions $d$
so can be discarded without affecting conformal invariance of the normal operator we are trying to construct. 
Vanishing of the variation of the Willmore invariant as given in Equation~\nn{slugmister} 
for general exact one-forms $\bar \Upsilon$ imposes the condition $\IIo=0$, or in other words umbilicity of the embedding.
This in turn implies that the
 Willmore invariant~$\bar B_3=0$.
Hence, because $\pi_1-1=0=\pi_2=\pi_3$  on singular Yamabe solutions, for umbilic embeddings $\Sigma^2\hookrightarrow (M^3,\cg)$  
we obtain 
\begin{align}\label{pi4ren}
\pi_4^{\rm ren}&:= 
 \hat n^a  \hat n^b \hat n^c  \nabla_{\hat n}
\nabla_a \nabla_b n_c|_\Sigma
+2H\!\hh P_{\hat n\hat n}
+3\hat n^a \hat n^b \nabla_{\hat n}\!\hh P_{ab }
\, .
\end{align}
This defines a conformal invariant of three dimensional singular Yamabe metrics for conformal manifolds whose boundary is
umbilically embedded. Moreover, since it is normal order three, it extracts the Neumann data.

\subsubsection{Expansion coefficients}
We can also compute the expansion coefficients for the singular Yamabe density $\csigma$ relative to a given CC metric $g_+$. Let us here  focus on $d=3$ dimensions. Let $g_{\rm GL}(\bar g)$ be the Graham--Lee compactified metric for a choice of $\bar g\in \bar \cg$ and $r$ the corresponding geodesic distance function $r$. Then 
$$
\sigma\stackrel{g_{\rm GL}}=
r\big(1+\tfrac12 H r -\tfrac{1}{3}[P_{\hat n\hat n}+K]r^2\big)
+\tfrac38 \bar B_3 r^4\log r
+o(r^4\log r)\, .
$$
When in addition the embedding is umbilic we have 
\begin{align}
\sigma
&\stackrel{g_{\rm GL}}=
r\big(1+\tfrac12 H r -\tfrac{1}{3}P_{\hat n\hat n}r^2\big)
\nonumber\\
&+\frac1{4!}\left[
3 
\sqrt{
\tfrac{\det 
g_{\sss\rm GL}^{\sss\rm SY}
}{\det g_{\sss\rm GL}}}\, 
\big(\hat n^a \hat n^b \hat n^c \nabla_a P_{bc}\big)^{g_{\sss\rm GL}^{\sss\rm SY}}
-3 \hat n^a \hat n^b \hat n^c \nabla_a P_{bc}
-2 H P_{\hat n\hat n}
\right]\hh r^4
\nonumber\\
&
\label{3Dexp}
+\frac1{5!}\left[4
\Big(
\tfrac{\det 
g_{\sss\rm GL}^{\sss\rm SY}
}{\det g_{\sss\rm GL}}
\Big)^{\frac23}
\Big(
 \hat n^a \hat n^b \hat n^c \hat n^d (\nabla_a \nabla_b P_{cd})
- (P_{\hat n\hat n})^2
\Big)^{g_{\sss\rm GL}^{\sss\rm SY}}
 \right.
 \\
 &
 \qquad\:\:\:
 -5H\hat n^a  \hat n^b \hat n^c   \nabla_{\hat n}\nabla_a\nabla_b n_c|_\Sigma
-10H^2  P_{\hat n\hat n}
-
4 \hat n^a \hat n^b \hat n^c \hat n^d (\nabla_a \nabla_b P_{cd})
\nonumber\\
&
\qquad\:\:\:
-20 H \hat n^a \hat n^b \hat n^c (\nabla_a P_{bc})
+4 (P_{\hat n\hat n})^2
+(\bar \nabla_a H)(\bar \nabla^a H)
 \Bigg]\: r^5
 + {\mathcal O}(r^6)\, ,\nonumber
\end{align}
where $g_{\sss\rm GL}^{\sss\rm SY}(\bar g)$ is the Graham--Lee compactified metric for the singular Yamabe metric~$g_+^{\rm SY}$. 
The above was obtained by comparing $\pi^{\rm ren}_4$ when computed with respect to~$g_{\sss\rm GL}^{\sss\rm SY}(\bar g)$  and~$g_+^{\rm SY}$. 
We also used that
 the mean curvature of a
 Graham--Lee compactified metric for the singular Yamabe metric vanishes.  
In addition we used the fifth order 
tensor of Section~\ref{beyond} to compute one term beyond the order of the Dirichlet-to-Neumann map.

\subsection{Fifth conformally compact boundary scalar}

Fifth CC boundary scalar curvatures  probe the failure of CC metrics to be singular Yamabe in dimensions~$d>4$. In four dimensions we can use them to study a generalized Willmore invariant as well as Neumann data for singular Yamabe solutions. Before writing down the requisite normal operator $\delta_5^{[4,5,7]}$, we define some conformally invariant operations.

\subsubsection{Conformally invariant operators and tensors}
First we need the fourth order Laplacian-squared 
type Paneitz operator of~\cite{Fradkin,Paneitz} intrinsic to $(\Sigma,\bar \cg)$ mapping $$\Gamma({\mathcal E}\Sigma[2-\tfrac{\bar d}2])\to \Gamma({\mathcal E}\Sigma[-2-\tfrac{\bar d}2])\, ,
$$
given for any
 $\bar g\in\bar \cg$ by
$$
\bar
{\sf P}_4 \stackrel {\bar g}{:=}\bar \Delta^2+\bar \nabla_a \circ (4 \bar P^{ab}-(\bar d-2) \bar J \bar g^{ab})\bar \nabla_b +\tfrac{\bar d -4}{2} \, \bar Q_{\bar d}\, ,
$$
where, for any dimension $\bar d$,
$$
\bar Q_{\bar d} \stackrel{\bar g}{:=} -\bar \Delta \bar J -2 \bar P^{ab}\bar P_{ab}
+ \tfrac{\bar d}2 \hh \bar J^2\, .
$$
The curvature $\bar Q_{4}$ is an example of a Branson $Q$-curvature~\cite{BO}.   In arbitrary dimensions~$\bar d$, we also need to define the operator
$$
\bar
{\sf P}_4^{(4)} \stackrel{\bar g}{:=}\bar \Delta^2\bar f+\bar \nabla_a \circ (4 \bar P^{ab}-2 \bar J \bar g^{ab})\bar \nabla_b  \, .
$$
If $\bar f$ is any function, then
\begin{multline}\label{Pvar}
\updelta^{\rm L} ({\sf P}_4^{(4)}\bar f)=(d\mm-\mm5)\big(
2\bar \nabla_{\bar \Upsilon}\mm \circ\mm \bar \Delta
+2 (\bar \nabla^a \bar \Upsilon^b)\bar \nabla_a\bar \nabla_b
\\
+(\bar \nabla^a \bar \nabla^b \bar \Upsilon_b) \bar \nabla_a
+2(d-1)\bar \Upsilon_a \bar P^{ab}\bar \nabla_b
\big)
 \bar f\, .
\end{multline}

We also need some second order differential operators. Firstly, 
 the ``partially massless'' operator $\bar {\sf M}^{(4)}$
 maps  
 $\Gamma(\odot^2_\circ T^*\Sigma[1])\to \Gamma(\odot^2_\circ T^*\Sigma[-1])$ according to (see~\cite{FradRep})
 \begin{equation}\label{PM}
\bar {\sf M}^{(4)} \mathring{\bar X}_{ab}\stackrel{\bar g}{ := }\bar \Delta   \mathring{\bar X}_{ab}
- \bar \nabla^c \bar \nabla_{(a} \mathring{\bar X}_{|c|b)_\circ}-\frac13
\bar \nabla_{(a} \bar \nabla^c \mathring{\bar X}_{|c|b)_\circ}\, .
\end{equation}
In dimensions $d>3$, we will need  a pair of second order, conformally invariant operators
$$\Delta^{\!\IIo}:\Gamma({\mathcal E}\Sigma[-1])\to \Gamma({\mathcal E}\Sigma[-4])\:\mbox{ and }
\Delta^{\!\!\IIIo}:\Gamma({\mathcal E}\Sigma[0])\to \Gamma({\mathcal E}\Sigma[-4])\, .
$$
These are defined for a choice of $g\in \cg$ by
\begin{align*}
\Delta^{\!\IIo}&\stackrel g{:=}
\, \bar \nabla_a\circ \IIo^{ab} \circ \bar \nabla_b 
-\tfrac{d-6}{d-2}  (\bar \nabla_a\IIo^{ab})\bar \nabla_b 
\\&\quad
-\tfrac 2{d-2} \big(\bar \nabla^a\bar \nabla^b \IIo_{ab}
+\tfrac{d-2}2\bar P_{ab}\IIo^{ab}
\big)
+\tfrac{2}{d-3} {\sf L}^{(3)}_{ab}\IIo^{ab}
\\
&=:\widehat{\bar \Delta}{}^{\!\IIo} +\tfrac{2}{d-3} {\sf L}^{(3)}_{ab}\IIo^{ab}
\,  ,\\[2mm]
\bar \Delta^{\!\!\IIIo}&\stackrel{g}{:=}
\bar \nabla_a\circ \IIIo^{ab} \circ \bar \nabla_b 
-\tfrac{d-5}{d-3} (\bar \nabla_a \IIIo^{ab})\bar \nabla_b\, .
\end{align*}
 In arbitrary dimensions~$d$ and any integer $\bar x\neq 2$, we  define the operator 
$$
\bar \Delta^{\!\!\IIIo}_{(\bar x)}\stackrel g{:=}
\bar \nabla_a\circ \IIIo^{ab} \circ \bar \nabla_b 
-\tfrac{\bar x-4}{\bar x-2} (\bar \nabla_a \IIIo^{ab})\bar \nabla_b\, .
$$
Also note that the linearized conformal variation of the operator  $(\bar \nabla_a \IIIo^{ab})\bar \nabla_b$, acting on a function $\bar f$, is given by
\begin{multline}\label{sausagedog}
\updelta^{\rm L}_{-4}\big((\bar \nabla_a \IIIo^{ab})\bar \nabla_b \bar f\big)=(d\mm-\mm3)\IIIo^{ab}\bar \Upsilon_a\bar \nabla_b \bar f
\\
=(d\mm-\mm 3)^2
\big(P^\top_{(ab)_\circ}\!\!-\mm\bar P_{(ab)_\circ}\!\! +\mm H\IIo_{ab}\big)\bar\Upsilon^a \bar \nabla^b \bar f
 \, .
\end{multline}

Along similar lines, in all dimensions $d>2$,   the second order, conformally invariant operator
$$
\bar \Delta^{\!\! K}:\Gamma({\mathcal E}\Sigma[0])\to \Gamma({\mathcal E}\Sigma[-4])\, ,
$$
 is defined for any $g\in \cg$ by
$$
\bar
\Delta^{\!\!K}\stackrel{g} {:=} K \bar \Delta +\tfrac{d-3}2 (\bar \nabla^a K)\bar \nabla_a\, .
$$
Also for any integer $\bar x$,  we call
$$
\bar \Delta^{\!\!K}_{(\bar x)}:= K \bar \Delta +\tfrac{\bar x-2}2 (\bar \nabla^a K)\bar \nabla_a\, .
$$

Note that, in dimensions $d>3$, the operator $K\bowtie\hat\delta_3$ of Equation~\nn{Kbowtie} can be written as the
manifestly invariant combination
$$
K\bowtie \hat \delta_3  \stackrel{g}= K \delta_3^{[2,3]}
-\Big(\tfrac{2(d-3)}{d-2} K^2 
-\tfrac{3}{d-3} \bar \Delta^{\!K} 
\Big) \delta_1 \, .
$$

In addition a pair of  bilinears are needed. These are conformal hypersurface invariants given by sections of ${\mathcal E}\Sigma[-4]$ that are 
defined for any $g\in \cg$ and dimension $d>2$, respectively,  by
\begin{align*}
\IIo\mm
\blacktriangleright\!\!
\blacktriangleleft \mm\IIo\:\:
&\stackrel{g}{:=}
\IIo^{ab} \bar \Delta \IIo_{ab}
+\tfrac{2(d-7)}{(d-2)} \IIo^{ab} \bar \nabla_a \bar \nabla^c \IIo_{cb}
\\&\:\:\:
-\tfrac{(d-12)(d-5)}{2(d-2)} (\bar\nabla^c \IIo^{ab})\bar \nabla_c \IIo_{ab}
+\tfrac{(d-5)(d-7)}{d-2} (\bar\nabla^c \IIo^{ab})(\bar \nabla_b \IIo_{ac})
\\&\:\:\:
+2(d-7)\bar P^{ab} \IIo_{ac} \IIo^c{}_b
-\bar J K\, ,
\\[1mm]
\bar \nabla.\IIo 
\blacktriangleright\!\!
\blacktriangleleft \bar \nabla.\IIo 
&
\stackrel{g}{:=}
(\bar\nabla_b \IIo^{ba})\bar \nabla^c \IIo_{ca}
-(d\!-\!2)\big(
(\bar\nabla^c \IIo^{ab})(\bar \nabla_c \IIo_{ab}
- \bar \nabla_b \IIo_{ac})
+\bar W^{abcd}\IIo_{ac}\IIo_{bd}\big)
 \, .
\end{align*}
We note the identity
$$
\bar \nabla.\IIo 
\blacktriangleright\!\!
\blacktriangleleft \bar \nabla.\IIo 
=-(d-2)\big(\bar W^{abcd}\IIo_{ac}\IIo_{bd}
+\tfrac12 \hh W^{\hat n abc} (W_{\hat n abc})^\top\big)\, ,
$$
so that in four dimensions
$
\bar \nabla.\IIo 
\blacktriangleright\!\!
\blacktriangleleft \bar \nabla.\IIo 
=-\tfrac12 \hh W^{\hat n abc} (W_{\hat n abc})^\top
$.
It also holds that
$$
\IIo
\blacktriangleright\!\!
\blacktriangleleft \IIo=
\left\{\begin{array}{cl}
-\frac32\Big( \bar{\sf L}^{ab} \big(\IIo_{c(a} \IIo^c{}_{b)\circ}\big)
-
\bar \nabla.\IIo 
\blacktriangleright\!\!
\blacktriangleleft \bar \nabla.\IIo
 \Big)\, ,
  & d=4\, ,\\[2mm]
  \IIo^{ab} \bar {\sf M}^{(4)} \IIo_{ab}
-\IIo^{ab}\bar W_{acbd}
\IIo^{cd}\, ,
  & d=5\, ,\\[2mm]
  \tfrac12 \bar \square^{(6)} K\, ,
   & d=7\, .\\[2mm]
\end{array}
\right.
$$
It follows in $d=4$ dimensions that
\begin{multline}\label{B4d}
\bar B_{4}\mm \stackrel{g}{:=} \mm\frac16\Big(
\bar {\sf L}^{ab}\mm \IIIo_{ab}
-\hh \IIo^{ab}\IVo{}_{ab}\mm
+\mm\IIIo^{ab} \IIo_{ac}\IIo^c{}_b
\\
+2\IIIo^{ab}\!\IIIo_{ab}
+\tfrac12K^2\mm
-\tfrac43
\IIo\mm
\blacktriangleright\!\!
\blacktriangleleft\mm \IIo
+
\bar \nabla.\IIo \mm
\blacktriangleright\!\!
\blacktriangleleft \mm\mm\bar \nabla.\IIo
\Big)\,\mm 
\end{multline}
agrees with the expression given in  Equation~\nn{B4}.
Importantly, we shall
define the tensor $\bar B_4$  in general dimensions $d$ by the
 above display, rather than  by Equation~\nn{B4}.
In particular, its linearized conformal variation is then given by 
\begin{align}\label{varB4}
\updelta^{\rm L}_{-4} \bar B_4= \tfrac{d-4}6\Big(
 \IIIo^{ab} \bar \nabla_a\bar \Upsilon_b
+2 (\bar \nabla^a \IIIo_{ab}) 
\bar \Upsilon^b
\Big)
=\tfrac{d-4}6 \,  \bar \Delta^{\!\!\IIIo}_{(4)} \bar\varpi\, .
\end{align}

\subsubsection{Fifth order boundary curvature}

We are now ready to build the fifth CC boundary curvature.
The operator 
$$\delta_5^{[4,5,7]}:\Gamma({\mathcal E}_0M[1])\to \Gamma({\mathcal E}M[-4])\, ,$$ is defined for $3<d\neq 4,5,7$ and some $g\in \cg$
by
\begin{align}\nonumber
 \delta^{\phantom {] I}\!\!\!\!}
 _{_{\scriptstyle 5}}
 {\!\!}^{[4,5,7]}
&\stackrel g{:=}\hat \delta_5
  - \tfrac{10}{d-7}\, \bar{\square}^{(6)}\circ \delta^{(7)}_3 
-\tfrac{1}{2(d-5)}\Big[
15\bar {\sf P}_4^{(4)}
-16\hh\IIo\mm\blacktriangleright\!\!
\blacktriangleleft \mm\IIo
-16 \hh\bar W^{abcd} \IIo_{ac}\IIo_{bd}
\Big]\circ\delta_1
\\&\:
+\tfrac{48}{d-4}\,
 \bar B_{4}
\delta_1
-\tfrac{15}{d-3}\Big[
\hh 
3\bar B_3 \delta_2
+\tfrac{2}{d-3} (\bar \nabla_a \IIIo^{ab})\bar \nabla_b \circ \delta_1
\Big]\, ,\label{polebit}
\end{align}
where
\begin{align*}
 \hat \delta_5&\stackrel{g}{:=}\hat n^a  \hat n^b \hat n^c \hat n^d 
 \hat n^e
 \nabla_a 
\nabla_b \nabla_c \nabla_d\nabla_e
+5 \big(H  \hat\delta_4 -H^2 \hat
 \delta_3
+{\mathcal S}\hh \delta_2\big)+(\hat {\mathcal T}+{\mathcal T}^{[1,2,3,4]})\hh\delta_1\, ,
\end{align*}
and
\begin{align*}
{\mathcal S}&:=
\hat n^a \hat n^b \hat n^c (\nabla_a P_{bc})
-2(\bar \nabla^a H)\bar \nabla_a
+\tfrac1{d-2} \Big[8(\bar\nabla_a \IIo^{ab})\bar \nabla_a
-3 (\bar \nabla^a \bar \nabla^b \IIo_{ab})\Big]
\, ,
\\[3mm]
\hat {\mathcal T} &:=
4 \hat n^a \hat n^b \hat n^c \hat n^d (\nabla_a \nabla_b P_{cd})
 +5 H \hat n^a \hat n^b \hat n^c (\nabla_a P_{bc})
-4 (P_{\hat n\hat n})^2-(\bar \nabla_a H)(\bar \nabla^a H)
\\&
+8 \bar P^{ab} \IIo_a{}^c \IIo_{cb}
-8\IIo^{ab}\hh\widehat{\IVo\hh\hh}{\!}_{ab}
+8K^2
\\&
+\Big[40 H \IIo^{ab}   \bar \nabla_a 
+5 (\bar \nabla^b P_{\hat n\hat n})
+30 \IIo^{ab} (\bar \nabla_a H) 
+ 20 H (\bar \nabla^b H)
-2 (\bar \nabla^b K)
\Big]\circ\bar \nabla_b
\\[2mm]
{\mathcal T}{}&{}^{[1,\ldots,4]}:=
-\tfrac{24}{d-4} \bar W^{abcd} \IIo_{ac}\IIo_{bd}
\\&
-\tfrac{1}{d-3}\Big[
30(\bar \nabla_a \IIIo^{ab})  \bar \nabla_b 
+8\bar \nabla^a \bar \nabla^b \IIIo_{ab}
+8\IIIo^{ab}\IIIo_{ab}
\Big]
\\&
+\tfrac{1}{d-2}\Big[
\tfrac{15}2 (\bar\nabla^a K)\bar \nabla_a
-60\IIo^{bc} (\bar \nabla^a \IIo_{ab})\bar \nabla_c
+40 \IIo^{bc}( \bar \nabla_b \bar \nabla^a \IIo_{ac})
-15 H (\bar \nabla^a \bar \nabla^b \IIo_{ab})
\\&\qquad\quad
+80(\bar \nabla^a \IIo^{bc})(\bar \nabla_a \IIo_{bc}-\bar \nabla_b \IIo_{ac})
-8 (\bar \nabla^a H)( \bar \nabla^b \IIo_{ab})
\\&\qquad \quad
+20 H^2 K -2 K^2
-8\bar \nabla.\IIo 
\blacktriangleright\!\!\blacktriangleleft \bar \nabla.\IIo
 + 80  \bar W^{abcd} \IIo_{ac}\IIo_{bd}
\Big]
\\&
+\tfrac{1}{(d-1)}
\Big[
\tfrac{15}2 (\bar\nabla^a K)\bar \nabla_a
\!+\!\tfrac{8}{3} \hh\IIo\!
\blacktriangleright\!\!
\blacktriangleleft \!\IIo
\!-\!32 \IIo^{bc} (\bar \nabla_b \bar \nabla^a \IIo_{ac})
\!-\!(\bar \nabla^a \IIo^{bc})(56\bar \nabla_a \IIo_{bc}\!-\!64\bar \nabla_b \IIo_{ac})
\\&\qquad\quad
+32 \bar P^{ab} \IIo_a{}^c \IIo_{cb}
-4K^2 -64 \bar W^{abcd} \IIo_{ac}\IIo_{bd}\Big]
\, .
\end{align*}
To construct the above expression we used that the operators
$$
K\hh \delta^{[3,2]}_3\, ,\quad
\IIIo^{ab}\IIo_{ab}\hh \delta_2\, ,\quad
\Delta^{\!\IIo}\circ\delta_2\, ,\quad
\bar \Delta^{\!K}\circ \delta_1\, ,\quad
\bar \Delta^{\!\mm\IIIo}\circ \delta_1\, ,
$$
are all invariant in dimensions $d\geq 4$, have transverse order $\leq 3$, and vanish when applied to a singular Yamabe scale. Hence they can be subtracted from any operator solving Problem~\ref{prob} and still yield a solution in the corresponding  critical dimensions. This is the case even when they appear as denominators of dimension-dependent poles.
A careful analysis of the remaining poles is given below, but already note that the residues of the poles at $d=3,4,5,7$ in Equation~\nn{polebit} are conformally invariant in respective dimensions $d=3,4,5,7$.  Also, the residues of the poles at $d=1,2,3,4$ in ${\mathcal T}^{[1,\cdots,4]}$ vanish in their respective dimensions.

Na\"ively one expects to obtain CC boundary scalars from the fifth order normal operator in any dimension $d\geq 5$. However, the presence of poles means that 
dimensions $d=5,7$ 
require special treatment.
But in dimensions 
$7\neq d\geq 6$ we already obtain a fifth CC boundary scalar
\begin{align*}
\pi_5 &\!:= \hat n^a  \hat n^b \hat n^c \hat n^d  \nabla_{\hat n} 
\nabla_a \nabla_b \nabla_c \nu_d|_\Sigma 
+5\big[H\pi_4-H^2 \pi_3+{\mathcal S}
\pi_2\big]+\tfrac{48}{d-4}\,
\bar B_{4}
\pi_1
\\&
-\tfrac{10}{d-7}\hh\bar\square^{(6)}\pi_3
+\tfrac{1}{2(d-5)}\Big[
60H\bar \square^{(4)}\pi_2
-
\Big(
15\bar {\sf P}_4^{(4)}
-16\hh\IIo\mm\blacktriangleright\!\!
\blacktriangleleft \mm\IIo
-16 \hh\bar W^{abcd} \IIo_{ac}\IIo_{bd}
\Big)\pi_1\Big]
\\&
-\tfrac{15}{d-3} \Big[3\bar B_3\pi_2
+\Big(
\tfrac12\bar \square^{(6)}\bar \Delta
+H^2\bar \Delta
-2H\IIo^{ab}\IIIo_{ab}
+3H\bar B_{3}
+\tfrac{2}{d-3} (\bar \nabla_a \IIIo^{ab})\bar \nabla_b
\Big) \pi_1\Big]
\\[1mm]&
+\tfrac{2}{d-2} \big[2\bar \square^{(6)}+5H^2\big](K\pi_1)
+\hat {\mathcal T}\pi_1+{\mathcal T}^{[1,\ldots,4]}\pi_1 
\, .
\end{align*}

For the case $d=7$, we use Equation~\nn{boxshift} to vary the tensor appearing in the numerator of the $d=7$ pole. This gives
$$-10\updelta^{\rm L}_{-2}(\bar\square^{(6)} \pi_3) = -10 (d-7) \bar \nabla_{\bar \Upsilon}\pi_3\, .
$$
Hence we condition on vanishing of $\pi_3$, which gives a fifth order  CC boundary curvature in $d=7$ dimensions, 
\begin{align*}
\pi_5|_{\pi_3=0,d=7}
&\!\stackrel{g}{:=}
 \hat n^a  \hat n^b \hat n^c \hat n^d  \nabla_{\hat n} 
\nabla_a \nabla_b \nabla_c \nu_d|_\Sigma 
+5\big[H\pi_4+\big(3H\bar \square^{(6)}+{\mathcal S}|_{d=7}\big)\, 
\pi_2\big]
\\&
-\tfrac{1}{4}\Big[
15\bar {\sf P}_4^{(4)}
-8(\bar \square^{(6)}K)
-16 \hh\bar W^{abcd} \IIo_{ac}\IIo_{bd}
\Big]\pi_1
+16\,
 \bar B_{4}
\pi_1
\\&
-\tfrac{15}{4} \Big[3\bar B_3\pi_2 
\!+\!\big(
\tfrac12
\bar \square^{(6)}\bar \Delta
+H^2\bar \Delta
\!-\!2H\IIo^{ab}\IIIo_{ab}
\!+\!3H\bar B_{3}
+\tfrac{1}{2} (\bar \nabla_a \IIIo^{ab})\bar \nabla_b
\big) \pi_1\Big]
\\[1mm]&
+\tfrac{2}{5} \big[2\bar \square^{(6)}+5H^2\big](K\pi_1)
+\big(\hat {\mathcal T}+{\mathcal T}^{[1,\ldots,4]}\big)|_{d=7}\:\pi_1 
\, .
\end{align*}

For the case $d=5$ we must examine the variation 
given in Equation~\nn{Pvar} of the operator ${\sf P}^{(4)}_4$ appearing in the numerator of the pole (note that the numerator proportional to $H\bar\square^{(4)} \pi_2$ cancels against the pole in the $5H\pi_4$ term).
It is clear that  for this to vanish for general conformal variations, we should condition upon vanishing of $\pi_1-1$. Hence for asymptotically hyperbolic structures in $d=5$ dimensions, we obtain the fifth order CC boundary curvature
\begin{align*}
\pi_5|_{\pi_1=1,d=5}
&\!\stackrel{g}{:=}\!
 \hat n^a  \hat n^b \hat n^c \hat n^d  \nabla_{\hat n} 
\nabla_a \nabla_b \nabla_c \nu_d|_\Sigma 
\!+\!5H\hat \pi_4
\!+\!5{\mathcal S}|_{d=5}\, 
\pi_2\!+\!\bar \square^{(6)} (5\pi_3\!-\!\tfrac43 K)
\!-\!5H^2 \pi_3
\\&
+\!48\bar B_{4}\mm
\! -\!\tfrac{45}2\bar B_{3}(\pi_2\mm+\mm H)
\!+\!15H\IIo^{ab}\IIIo_{ab}
\!+\!\tfrac{10}{3} H^2K
\!+\!(\hat {\mathcal T}\hh 1)|_{d=5}\mm+\mm{\mathcal T}^{[1,\ldots,4]}|_{d=5} 
\, ,
\end{align*}
where
\begin{multline*}
\hat \pi_4=
 \hat n^a  \hat n^b \hat n^c \hat n^d \nabla_a 
\nabla_b \nabla_c \nu_d|_\Sigma
+2H \pi_3
+ P_{\hat n\hat n} \pi_2
-3\IIo^{ab}\IIIo_{ab}+\tfrac94\bar B_{3}
\\
-\tfrac{4}{3} \!\hh H K
 +(\bar\nabla^a\bar \nabla^b \IIo_{ab})
+2H\!\hh P_{\hat n\hat n}
+3\hat n^a \hat n^b (\nabla_{\hat n}\!\hh P_{ab })\, .
 \end{multline*}

\subsubsection{Four dimensional Dirichlet-to-Neumann map}

We now study the $d=4$ case. The residue of the $d-4$ pole in $\pi_5$ recovers the obstruction density $\boldsymbol{B}_4$. We must now analyze the failure of the numerator $48\bar B_4 \pi$ to be invariant in general dimensions $d$. While this is necessarily proportional to $(d-4)$, there is the freedom to add terms proportional to $d-4$ to the numerator to minimize the failure of invariance. This is why we chose to define $\bar B_4$ in dimensions not equal to $4$ by Equation~\nn{B4d} rather than~\nn{B4}.
The superumbilic condition requires that both $\bar B_4$ (see Equation~\nn{B4d}) and the variation of its extension to $d$ dimensions (see Equation~\nn{varB4}) vanish. For this, it suffices to assume that the embedding is both umbilic and Fialkow flat. 
We can then remove the $d=4$ pole, and obtain the Dirichlet-to-Neumann map in  $d=4$ dimensions and for umbilic, Fialkow-flat embeddings,
\begin{multline}\label{pi5ren}
\pi_5^{\rm ren}=
\hat n^a  \hat n^b \hat n^c \hat n^d  \nabla_{\hat n}\nabla_a\nabla_b\nabla_c n_d|_\Sigma
+
4 \hat n^a \hat n^b \hat n^c \hat n^d (\nabla_a \nabla_b P_{cd})
\\
 +5 H \hat n^a \hat n^b \hat n^c (\nabla_a P_{bc})
-4 (P_{\hat n\hat n})^2
-(\bar \nabla_a H)(\bar \nabla^a H)
\, .
\end{multline}
Let us remark that  the umbilic and  Fialkow-flat 
conditions in four bulk dimensions ensure that the metric $g_+^{\rm SY}$ is asymptotically Poincar\'e--Einstein, meaning
that it solves the Poincar\'e--Einstein to the highest possible order determined by formal asymptotics alone.

\subsubsection{Beyond the Dirichlet-to-Neumann map}
\label{beyond}

In dimension $d=3$ we can also use the fifth order operator~\nn{polebit} to probe higher order terms in the expansion~\nn{3Dexp}. In order that there are no log terms, we must again assume that the embedding is umbilic. This also ensures the Willmore invariant $\bar B_3=0$. To discard the $d=3$ pole in Equation~\nn{polebit}, we  must examine the conformal variation of 
the operator $3\bar B_3 \delta_2
-\tfrac{2}{d-3} (\bar \nabla_a \IIIo^{ab})\bar \nabla_b \circ \delta_1$, to which end one can consult Equations~
\nn{slugmister} and~\nn{sausagedog}. Since $\IIIo$ is identically zero in $d=3$ dimensions, and we already have imposed the umbilic condition $\IIo=0$, it follows that even upon removing the $d=3$ pole and acting on the singular Yamabe defining density, we have that
\begin{multline*}
\pi_5^{\rm ren}|_{d=3}\stackrel{g}{:=}
\hat n^a  \hat n^b \hat n^c \hat n^d  \nabla_{\hat n}\nabla_a\nabla_b\nabla_c n_d|_\Sigma
+5H\hat n^a  \hat n^b \hat n^c   \nabla_{\hat n}\nabla_a\nabla_b n_c|_\Sigma
+10H^2\hat n^a  \hat n^b  P_{ab}\\
+
4 \hat n^a \hat n^b \hat n^c \hat n^d (\nabla_a \nabla_b P_{cd})
+20 H \hat n^a \hat n^b \hat n^c (\nabla_a P_{bc})
-4 (P_{\hat n\hat n})^2
-(\bar \nabla_a H)(\bar \nabla^a H)
\end{multline*}
is  conformally invariant. The above fifth order tensor was used to compute the fifth order term in Equation~\nn{3Dexp}.

\subsubsection{Expansion coefficients}\label{EC5}

We can now compute further expansion coefficients for the singular Yamabe density relative to a given CC metric. Firstly let us focus on $d=4$ dimensions. As per the discussion in Section~\ref{EC}, in the case that the obstruction density~$\bar{ \cB}_4$ is non-vanishing, we can determine the first four expansion coefficients and the coefficient of the first log term. Given $\bar g\in \bar \cg$, we find
\begin{multline*}
\sigma\stackrel{g_{\rm GL}}=
r\big(1+\tfrac12 H r -\tfrac{1}{3}[P_{\hat n\hat n}+\tfrac12 K]r^2
\\
-\tfrac1{24}\big[
2H\!\hh (P_{\hat n\hat n}\!-\!K)
+3\hat n^a \hat n^b (\nabla_{\hat n}\!\hh P_{ab })
+\tfrac32 (\bar \nabla^a \bar \nabla^b \IIo_{ab})
-6\IIo^{ab}\IIIo_{ab}+9\bar B_{3}
\big]r^3\big)
\\
+\tfrac25 \bar B_4 r^5\log r
+o(r^5\log r)\, .
\end{multline*}

When in addition the embedding is superumbilic we can give a formula for the first five expansion coefficients
\begin{align*}
\sigma
&\stackrel{g_{\rm GL}}=
r\Big(1+\tfrac12 H r -\tfrac{1}{3}P_{\hat n\hat n}r^2
-\tfrac1{24}\big[
2H\!\hh P_{\hat n\hat n}
+3\hat n^a \hat n^b (\nabla_{\hat n}\!\hh P_{ab })
\big]r^3\Big)
\\&
+\frac1{5!}\left[
4\Big(
\tfrac{\det 
g_{\sss\rm GL}^{\sss\rm SY}
}{\det g_{\sss\rm GL}}
\Big)^{\frac23}
\Big(
 \hat n^a \hat n^b \hat n^c \hat n^d (\nabla_a \nabla_b P_{cd})
- (P_{\hat n\hat n})^2
\Big)^{g_{\sss\rm GL}^{\sss\rm SY}}
 \right.
 -4 \hat n^a \hat n^b \hat n^c \hat n^d (\nabla_a \nabla_b P_{cd})
 \\&\quad\quad\;\;
-5 H \hat n^a \hat n^b \hat n^c (\nabla_a P_{bc})
+4 (P_{\hat n\hat n})^2
+(\bar \nabla_a H)(\bar \nabla^a H)
 \Bigg]\: r^5
 + {\mathcal O}(r^6)\, .
\end{align*}

We can also give expansion formul\ae\ in higher dimensions. 
For simplicity we focus on $d=6,8,9,\ldots$, although only slight adjustments are needed to handle the cases $d=5,7$. 
As discussed in Section~\ref{EC}, we can then use vanishing of $\pi_1-1$, $\pi_2,\ldots, \pi_5$ to determine coefficients $s_0,\ldots, s_5$. For a choice of $\bar g\in \bar \cg$ this gives
\begin{multline*}
\hspace{-.3cm}\sigma\stackrel{g_{\sss\rm GL}}= r\Big\{1
+\tfrac12 Hr
-\tfrac13\big[P_{\hat n\hat n}+\tfrac1{d-2} K\big]r^2
\\
-\tfrac1{24}\Big[
3\hat n^a \hat n^b (\nabla_{\hat n}\!\hh P_{ab })
\!+\!2H\!\hh P_{\hat n\hat n}
\!-\!\tfrac3{d-3} \big(2\IIo^{ab}\IIIo_{ab}\!-3\bar B_{3}\big)
\!+\!\tfrac1{d-2}(3\bar\nabla^a\bar \nabla^b\IIo_{ab}
\!-\!4HK)
\Big]r^3
\\
-\tfrac1{120}
\Big[4 \hat n^a \hat n^b \hat n^c \hat n^d (\nabla_a \nabla_b P_{cd})
+ 5 H \hat n^a \hat n^b \hat n^c (\nabla_a P_{bc})
-4 (P_{\hat n\hat n})^2
\phantom{operpiugougyaootoo}\hh
\\[-1mm]
-(\bar \nabla_a H)(\bar \nabla^a H)
+8 \bar P^{ab} \IIo_a{}^c \IIo_{cb}
-8\IIo^{ab}\hh\widehat{\IVo\hh\hh}{\!}_{ab}
+8K^2\phantom{splattotozz}
\\[2mm]
\!\!\!\!
+\tfrac{8}{d-5}\big(
\IIo\mm\blacktriangleright\!\!
\blacktriangleleft \mm\IIo
+ \bar W^{abcd} \IIo_{ac}\IIo_{bd}
\big)
+\tfrac{48}{d-4}\,
( \bar B_{4}-\tfrac12  \bar W^{abcd} \IIo_{ac}\IIo_{bd})
\phantom{suck}
\\
\!\!\!\!
 -\tfrac{1}{d-3} \big(
 8\bar \nabla^a \bar \nabla^b \IIIo_{ab}
+8\IIIo^{ab}\IIIo_{ab}
-30H\IIo^{ab}\IIIo_{ab}
-45H\bar B_{3}
\big) 
\phantom{spurt0}
\\\!\!\!\!\!\!\!\!\!\!\!
+\tfrac{1}{d-2} \Big(4 \bar\square^{(6)}K+10H^2K
+40 \IIo^{bc}( \bar \nabla_b \bar \nabla^a \IIo_{ac})
-15 H (\bar \nabla^a \bar \nabla^b \IIo_{ab})\!\!\!
\\
+80(\bar \nabla^a \IIo^{bc})(\bar \nabla_a \IIo_{bc}-\bar \nabla_b \IIo_{ac})
-8 (\bar \nabla^a H)( \bar \nabla^b \IIo_{ab})
+20 H^2 K\hspace{-1.6cm}
\\
 -2 K^2
-8\bar \nabla.\IIo 
\blacktriangleright\!\!\blacktriangleleft \bar \nabla.\IIo
 + 80  \bar W^{abcd} \IIo_{ac}\IIo_{bd}
\Big)
\phantom{groggo}
\\
\hspace{-2.7cm}
+\tfrac{1}{(d-1)}
\Big(
\tfrac{8}{3} \hh\IIo\!
\blacktriangleright\!\!
\blacktriangleleft \!\IIo
-32 \IIo^{bc} (\bar \nabla_b \bar \nabla^a \IIo_{ac})
\phantom{pumpkinnos}
\\
-(\bar \nabla^a \IIo^{bc})(56\bar \nabla_a \IIo_{bc}
-64\bar \nabla_b \IIo_{ac})
+32 \bar P^{ab} \IIo_a{}^c \IIo_{cb}
\\
-4K^2 -64 \bar W^{abcd} \IIo_{ac}\IIo_{bd}\Big)
\Big]
r^4
\Big\}
+{\mathcal O}(r^6)\, .\phantom{applewinebarrelleakmi}
\end{multline*}


\section*{Acknowledgements}
We thank John Chae for a collaboration during early stages of this work.
A.R.G. and A.W. 
 acknowledge support from the Royal Society of New Zealand via Marsden Grant  19-UOA-008. 
A.W.~was also supported by  Simons Foundation Collaboration Grant for Mathematicians ID 686131. J.K. and A.W. thank the Banff International Research Station for providing a stimulating research environment for part of this work.

\end{document}